\documentclass[11pt]{amsart}
\usepackage[T1]{fontenc}
\usepackage{amssymb, amsthm, amsmath}

\usepackage[mathscr]{eucal}

\newtheorem{theorem}{Theorem}[section]
\usepackage{graphicx}%
\usepackage{hyperref}
\newtheorem{introtheorem}{Theorem} 

\newtheorem{lemma}[theorem]{Lemma}
\newtheorem{fact}[theorem]{Fact}
\newtheorem{cor}[theorem]{Corollary}
\newtheorem{claim}[theorem]{Claim}
\newtheorem{proposition}[theorem]{Proposition}

\theoremstyle{definition}

\newtheorem{rmk}[theorem]{Remark}
\newtheorem{defnn}[theorem]{Definition}
\newcommand\vlabel{\label}
\newcommand{\la}{\langle}
\newcommand{\ra}{\rangle}

\newcommand{\CY}{{\mathcal Y}}
\newcommand{\Cq}{{\mbox{$C_{q_0}$}}}

\newcommand{\CL}{{\mathcal L}}

\newcommand{\CN}{{\mathcal N}}

\newcommand{\CR}{{\mathcal R}}
\newcommand{\CM}{{\mathcal M}}

\newcommand{\CF}{{\mathcal F}}

\newcommand{\CS}{\mathcal S}
\newcommand{\sub}{\subseteq}
\newcommand{\pf}{\proof }

\newcommand{\x}{\bar x}
\newcommand{\y}{\bar y}
\newcommand{\0}{\emptyset}

\newcommand{\fleq}{\preccurlyeq}

\newcommand{\C}{\mathcal C}

\newcommand{\CX}{{\mathcal X}}
\newcommand{\acl}{\operatorname{acl}}
\newcommand{\dcl}{\operatorname{dcl}}
\newcommand{\Macl}{\operatorname{acl_{\CM}}}
\newcommand{\Nacl}{\operatorname{acl_{\CN}}}
\newcommand{\Mdcl}{\operatorname{dcl_{\CM}}}

\newcommand{\mr}{\operatorname{MR}}

\newcommand{\st}{\operatorname{st}}
\newcommand{\cl}{\operatorname{cl}}
\renewcommand{\leq}{\leqslant}
\renewcommand{\geq}{\geqslant}
\newcommand{\reals}{\mathbb R}
\newcommand{\Mdim}{\dim_{\CM}}
\newcommand{\Ndim}{\dim_{\CN}}
\newcommand{\eps}{\epsilon}

\newlength{\intlength}

\usepackage{endnotes}

\title{One dimensional structures in o-minimal theories}

\author[A. Hasson]{Assaf Hasson$^*$}
\thanks{$^*$Supported by the EPSRC grant no. EP C52800X 1}
\address{University of Oxford, Mathematical Institute, 24-29 St Giles', Oxford, OX1 3LB, UK}
\email{hasson@maths.ox.ac.uk }

\author[A. Onshuus]{Alf Onshuus}
\address{Universidad de los Andes,
Departemento de Matem\'aticas, Cra. 1 No 18A-10, Bogot\'{a},
Colombia} \email{onshuus@gmail.com}

\author[Peterzil]{Ya'acov Peterzil}
\address{Department of Mathematics,
University of Haifa, Haifa, ISRAEL} \email{kobi@math.haifa.ac.il}

\begin{document}

\begin{abstract} We prove the Zil'ber Trichotomy Principle for
all 1-dimensional structures which are definable in o-minimal ones.

On the level of types, we prove the following: Let $\CM$ be an
o-minimal structure and let $\CN$ be any structure definable in
$\CM$. Assume that $p$ is an $\CN$-type which is 1-dimensional in
the sense of $\CM$. Then, essentially, $\CN$ induces on $p$ either a
trivial $acl_{\CN}$-geometry, or the structure of a pure (possibly
ordered) $\CN$-definable vector space over a division ring, or an
$\CN$-definable o-minimal expansion of a real closed field.

In particular, if $\CN$ is stable and 1-dimensional then the $acl$
geometry is necessarily locally modular.

Along the way, we develop a fine intersection theory for definable
curves in o-minimal structures. 
\end{abstract}

\maketitle

The present work is a first step in an attempt to classify
structures interpretable in o-minimal theories. Such structures may
be stable (e.g., algebraically closed fields of characteristic zero,
compact complex manifolds) or unstable (e.g., expansions of real
closed fields, ordered vectors spaces). The ultimate goal would be
to treat in some uniform manner all of these examples by exploiting
the ambient  o-minimal environment.

\begin{defnn}
A structure $\CN$ is \emph{definable} in an o-minimal structure
$\CM$ if it is interpretable in the real sorts of $\CM$ (i.e. not in
$\CM^{eq}$). More precisely, the universe of $\CN$, call it $N$, is
an $\CM$-definable set and its atomic relations are certain
$\CM$-definable subsets of $N^k$, $k\in \mathbb N$. We say that
$\CN$ is \emph{definable of dimension $k$ in $\CM$} if $k$ is the
minimal possible dimension of such an $\CM$-definable set
$N$.
\end{defnn}

If $\CN$ is definable  in an o-minimal structure $\CM$, quite a few
good properties are induced on $\CN$. Here is a partial list:
\begin{enumerate}
\item$\CN$ has the non-independence property, NIP.
\item $\CN$ is super-rosy of finite \th-rank, \cite{On}
\item There is a bound on uniformly definable families of
finite sets in $\CM$ and therefore also in $\CN$. Said differently,
 $\CN$ eliminates the
$\exists^{\infty}$ quantifier.
\item If $N$ is a $\CN$-minimal (i.e. $N$ has no
$\CN$-definable infinite subsets of smaller $\CM$-dimension) then $acl_{\CN}$
(the algebraic closure in the sense of $\CN$) satisfies
the exchange property (see \cite{PePiSt2}). Combined with (3), this means that
$\CN$ is a geometric structure (see \cite{HrPi1} for more
 on geometric structures).
\item  If $\CN$ is stable then it has finite $U$-rank,
 with $U(\CN)\leq \dim_{\CM}(N)$,
\cite{Ga}.
\end{enumerate}

Let us consider the following setting:
{\em Fix $\CM$,  a sufficiently saturated o-minimal structure
 with a dense underlying order}, and assume that
 $\CN$ is a geometric structure definable
in $\CM$.

Our main interest  is to extract algebraic
information from the geometric structure of $\CN$. The idea is to 
follow Zilber's division of geometric structures $\CN$ into
three types:

\begin{description}
    \item[Degenerate] $\Nacl(A)=\bigcup\{\Nacl(a): a\in A\}$ for
    all $A\subset \mathcal N$.
    \item [Linear] Every normal definable family of plane curves
    in $\mathcal N$ is 1-dimensional (in the sense of $\CN$), 
    but $\mathcal N$ is
    non-degenerate.
    \item[Rich] There exists a normal 2-dimensional (in 
the sense of $\CN$) definable family
 of plane curves.
\end{description}

Clearly, degenerate geometries do not allow the existence of
infinite definable groups, and indeed by compactness not even the
existence of type definable ones. Therefore there is no point
looking for algebraic data in degenerate structures. But the
Linear/Rich dichotomy has been subject to much research. The
underlying thesis, sometimes 
known as Zilber's Principle, is that linear geometries
should arise from pure linear structures (e.g a vector space with no
additional structure) and that in contexts of "topological flavour"
rich combinatorial geometries arise from the geometry of interpretable
fields.

In the stable case, though not in general true (see \cite{HR1}),
 Zilber's Principle
has been first proved in \cite{HZ} for strongly minimal Zariski geometries,
 and
extended later to several related contexts. In the o-minimal world a local
 version
 of the principle was proved in \cite{PeSt2}.

A natural conjecture to make  in this context is:
\\

\noindent{\bf Conjecture } 
{\em If $\CN$ is a geometric structure definable in 
an o-minimal $\CM$  then Zilber's Principle holds for $\CN$.
In particular, if $\CN$ is rich then it interprets a field.}
\\

Note that the above conjectre, if true, will cover some remaining
open cases of Zilber's original conjecture,  for strongly minimal
structures interpretable in algebraically closed fields of charateristic zero
(See Rabinovich's work \cite{Ra} on the main cases of this question).

Our present work can be seen as a proof of the above conjecture
in the case where $\CN$ is definable in $\CM$ and
$\dim_{\CM}(N)=1$ (see \cite{HaKo} and \cite{PeSt} for the investigation
of certain  cases of definable, stable, two-dimensional structures). 
 First note that if 
 $\CN$ is  1-dimensional and definable in $\CM$
then  every $\CN$-definable subset of $N$ is
 either finite or has
dimension one. Thus $N$ itself is an $\CN$-minimal set and therefore
a geometric structure.
The main bulk of our work goes into
proving that if  $\CN$ is one-dimensional, definable in 
$\CM$, and stable then  it cannot be
rich. Let us see, in view of  Zilber's Principle above, why we should
indeed expect this to be true:

 Assume that $\CN$ is rich, then 
a field $K$ should be interpretable in $\CN$ (if we 
believe the above principle). Because $\CN$ is superstable
the field $K$ must be algrebraically closed, and because it is interpretable in an o-minimal theory it has to be of characteristic zero. 
Moreover, if $\CM$ expands a real closed field $R$
then $K$ is definably isomorphic to the algebraic closure of $R$ 
and in particular its o-minimal dimension is two 
(for both statements see \cite{OtPePi}). It follows (\cite{PeSt}) that $K$, 
with all its $\CN$-induced structure, is a pure field and therefore 
strongly minimal. However, because $K$ is interpretable in 
the geometric structure $\CN$, its strong minimality implies
that $\dim_{\CM}(K)=\dim_{\CM}(N)=1$, a contradiction. 
(The assumption that $\CM$ expands a real closed
field can actaully be avoided).

The same argument shows, assuming that
the conjecture is true: If
$\CN$ is a rich strongly minimal structure which is interpretable in an 
o-minimal
one then necessarily $\dim_{\CM}(N)=2$. Interestingly,
the remaining open cases in proving Zilber's conejcture for strongly minimal
structures  in algebraically closed fields are exactly
those in which the underlying universe of $\CN$ is of 
{\em algebraic dimension}
bigger than one.

Returning to our present work, after proving that a stable 
one-dimensional  structure definable in $\CM$ cannot be rich,
the analysis of the linear stable case reduces to
classical structure theorems on 1-based theories of $U$-rank 1. 
In the unstable case, we
rely on results from \cite{HaOn} and \cite{PeSt2} to finish our analysis.
Taken together, our main theorem is:

\begin{introtheorem}\vlabel{main} Let $\CN$ be a 1-dimensional structure
 definable in
an o-minimal one. Then one of the following hold:

\begin{enumerate}
    \item $\CN$ is degenerate.
    \item $\CN$ is linear: There exists $X\subseteq N^{eq}$ such that either $X$ is
strongly minimal and locally modular (whose $\CN$-definable structure arises from a
definable vector space) or $X$ is an o-minimal group interval, whose $\CN$-definable structure arises from an ordered vector space.
    \item $\CN$ is rich and interprets a real closed field.
\end{enumerate}
\end{introtheorem}

Because of the flexibility in patching together any two o-minimal structures into a new one, it is clear that on the global level we cannot hope for sharper results. To get the more informative statement,  our theorem will be formulated in terms of certain types in an arbitrary $\CM$-definable structure $\CN$.  
We say that an $\CN$-type $p$ is 
{\em $\CM$-1-dimensional} if it is contained in
some $\CN$-definable set whose $\CM$-dimension is one.
Because $\CN$ has finite \th-rank,
 a natural first step in any analysis
 of such an $\CN$ would be a classification 
of its \th-rank 1 types. Here  we complete this
 classification for the case when \th-rank minimal types
are $\CM$-1-dimensional. The work of \cite{HaOn} suggests
 that for unstable types this classification is not far from being 
   sufficient.  For a non-algebraic type
$p\in S_1(N)$ say that:
\begin{description}
\item[Trivial] $p$ is trivial if for every finite set of parameters
$A\subseteq N$, and $a_1,\ldots, a_n,b\models p$, if $b\in
\Nacl(Aa_1...a_n)$ then $b\in \bigcup_{i=1}^n \Nacl(Aa_i)$.
\item[Rich] $p$ is rich if there exist a finite set $A$, an element $b\models p$ with
$\dim_{\CN}(b/AN)=1$ and an (almost) normal family
 $\CF$ of plane curves, $\CN$-definable over the set $NA$, such
  that $\{f\in \CF:\la b,b \ra\in f\}$ is infinite.
\item[Linear] $p$ is linear if $p$ is not trivial and not rich.
\end{description}


 Note that the notion of
richness is definable
 in the following sense: If $p$ is rich, witnessed by $b$ and a family
 of curves $\CF_{\bar a}$
 (over parameters in $\bar a$) then, because $\Ndim(b/Na)=1$,
 there exists an $\bar a$-definable infinite set $X_{\bar a}$ containing
 $ b$ such that whenever $b'\in X_{\bar a}$ the pair
 $\la b',b'\ra$ belongs to infinitely many curves from $\CF_{\bar a}$.
  We can now find a formula $\theta\in p$ such that for all
$b'\models \theta$ there exist a tuple $\bar a'$, $\CF_{\bar a'}$
and $X_{\bar a'}$ with the same properties. We  will say that
$\theta$ isolates the richness of $p$.

The local version of Theorem \ref{main} is then:

\begin{introtheorem}\vlabel{mainlocal}
Assume that $\CN$ is a definable structure in an o-minimal $\CM$
and that $p$ is a complete 1-$\CM$-dimensional $\CN$-type over a model
$\CN_0\sub \CN$. Then:
\begin{enumerate}
\item $p$ is trivial (with respect to $acl_{\CN}$).

Or, there exists an $\CN$-definable equivalence relation $E$ with
finite classes such that one of the following holds:

\item $p$ is linear, in which case either
\begin{enumerate}
 \item $p/E$ is a generic type of a strongly minimal $\CN$-definable
one-dimensional group $G$, and the structure which $\CN$ induces on
$G$ is locally modular.  In particular, $p$ is strongly minimal. Or,
 \item for every $a\models p$ there exists
 an $\CN$-definable ordered group-interval $I$ containing $a/E$.
 The structure which $\CN$ induces on $I$ is a reduct of an ordered
vector space over an ordered division ring.
\end{enumerate}
\item $p$ is rich: For every $a\models p$ there exists
 an $\CN$-definable real closed field $R$ containing $a$
 and the structure which $\CN$ induces on $R$ is o-minimal.
\end{enumerate}
\end{introtheorem}

The above statement gives a more precise meaning to Theorem \ref{main} and clearly strengthen it.

  The results
of this paper are true for $\CN$ \emph{definable} in arbitrary dense
o-minimal structures (and by elimination of imaginaries, we may
replace "definable" with "interpretable", if $\CM$ expands an
o-minimal group). However, in order  to keep the exposition cleaner,
throughout the main part of the work, we will assume that in fact
$\CM$ expands an o-minimal field. In Appendix B we prove the results
to arbitrary o-minimal structures.

\noindent {\bf The structure of the paper } Section 1 is quite technical. In it
we develop a fine theory of tangency and transversality for curves
in definable families. An important ingredient in the whole argument
turns out to be the theory of limit sets as developed by v.d. Dries
in \cite{vdDr2} and we review it here and in Appendix A. In Section
2 we prove that every stable definable 1-dimensional structure is
necessarily 1-based. In Section 3 we prove the main theorem, for
structures and for types. As mentioned above, in Appendix B we
generalise the result from expansions of real closed fields to
arbitrary o-minimal structures.

\section{Some elements of intersection theory}\label{Intersection}

In this section we develop elements of intersection theory which
will be used in the proof of Theorem \ref{main}. Some parts of this
theory were already developed in \cite{PeSt2} and we will use
several results and proofs from that paper. However, the setting
there allowed us to ``trim'' the original family of curves as the
proof progressed, by working locally in regions where the original
family is ``well-behaved''. In the current setting this is
impossible (because the ordering is not assumed to be definable in
$\CN$) and therefore one needs to develop a finer intersection
theory. We mainly concentrate on counting the number of intersection
points of plane curves.

Although the treatment of the strongly minimal case will focus on
reaching a contradiction (from the assumption that the structure is
one dimensional and locally modular) the results proved in this
section are true in any o-minimal structure with a family of curves
as given below. Some of the results are stated under the assumption
that the o-minimal structure expands a field. In Appendix B we will
show how to avoid this assumption.

\subsection{Intersections of two curves}

We use the term {\em curve} to denote any $\CM$-definable one-dimensional set.
 A {\em plane curve} is a curve in $M^2$.

\begin{defnn}
{\em Given an $\CM$-definable plane curve and $p\in X$, we say that
{\em $X$ is $C^0$ at $p$} if $X$ is locally $\CM$-homeomorphic to an
open interval (we also say that $p$ is a $C^0$-point of $X$).}
\end{defnn}

 It is not difficult to see that in this case
$X$ divides $M^2$, sufficiently close to $p$,  into two definably
connected components, $W_1, W_2$.  Namely, for every sufficiently
small rectangular neighbourhood  $W$ of $p$, the set $W\setminus X$
has two definably connected components, one having the same germ as
$W_1$ at $p$ and the other as $W_2$.

The following definitions come from 2.13 of \cite{PeSt2}:

\begin{defnn}{\em
If $X$ and $Y$ are two $\CM$-definable plane curves and $p$ is a
$C^0$-point on both, we say that $X$ and $Y$  {\em touch each other at
  $p$} (see Figure 1 below) if there is some rectangular
neighbourhood $U$ of $p$ such that $U\setminus X$ has
two components $W_1,W_2$ and either $Y\cap U\cap W_1=\emptyset$ or
$Y\cap U\cap W_2=\emptyset$.}

\end{defnn}

\begin{center}%
\begin{table}[h]
\begin{tabular}{ll}%
{\parbox[b]{2in}{\begin{center}
\includegraphics[height=1.2in,width=2in]
{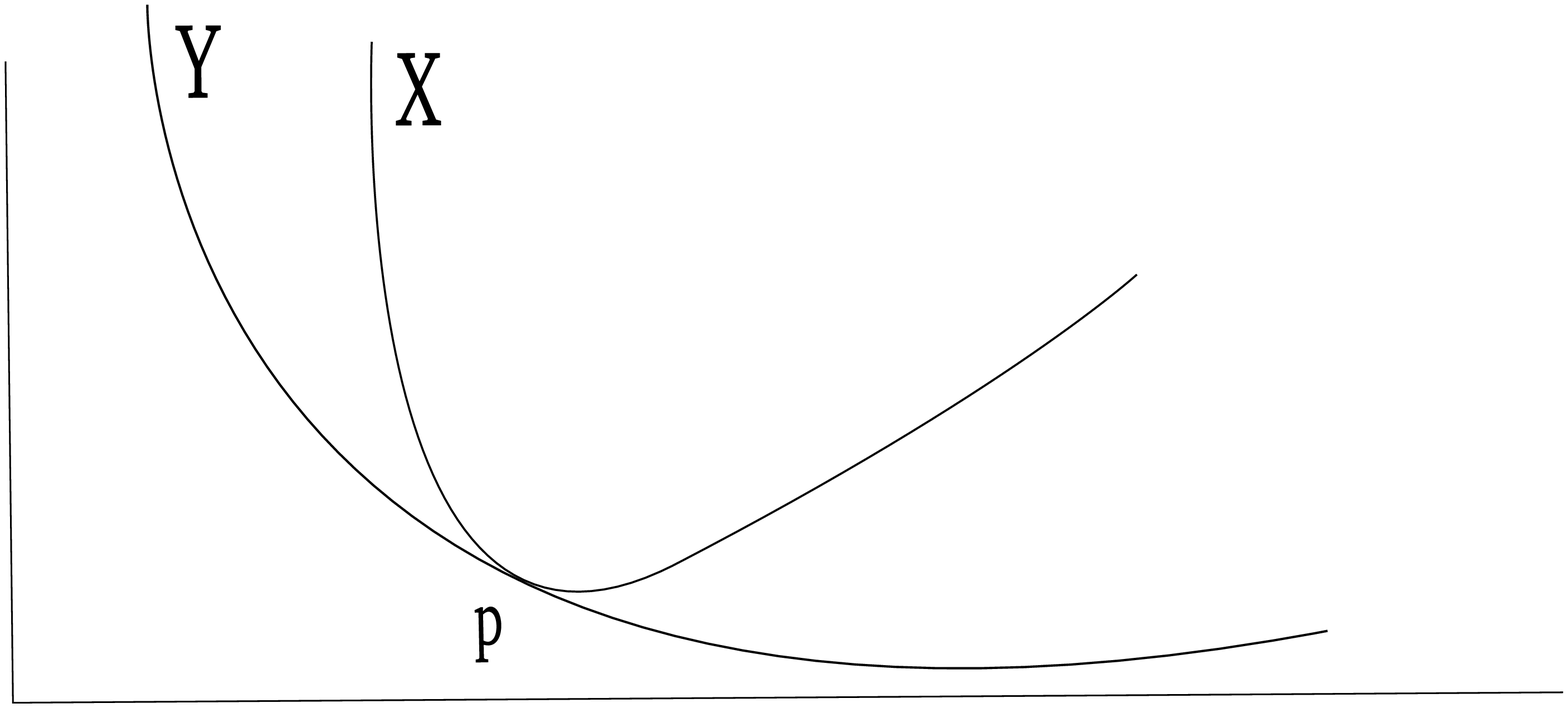}\\
Figure 1: \\ $X$ touches $Y$ at $p$\end{center}}}%
{\parbox[b]{4in}{
\begin{center}
\includegraphics[height=1.2in,width=2in]{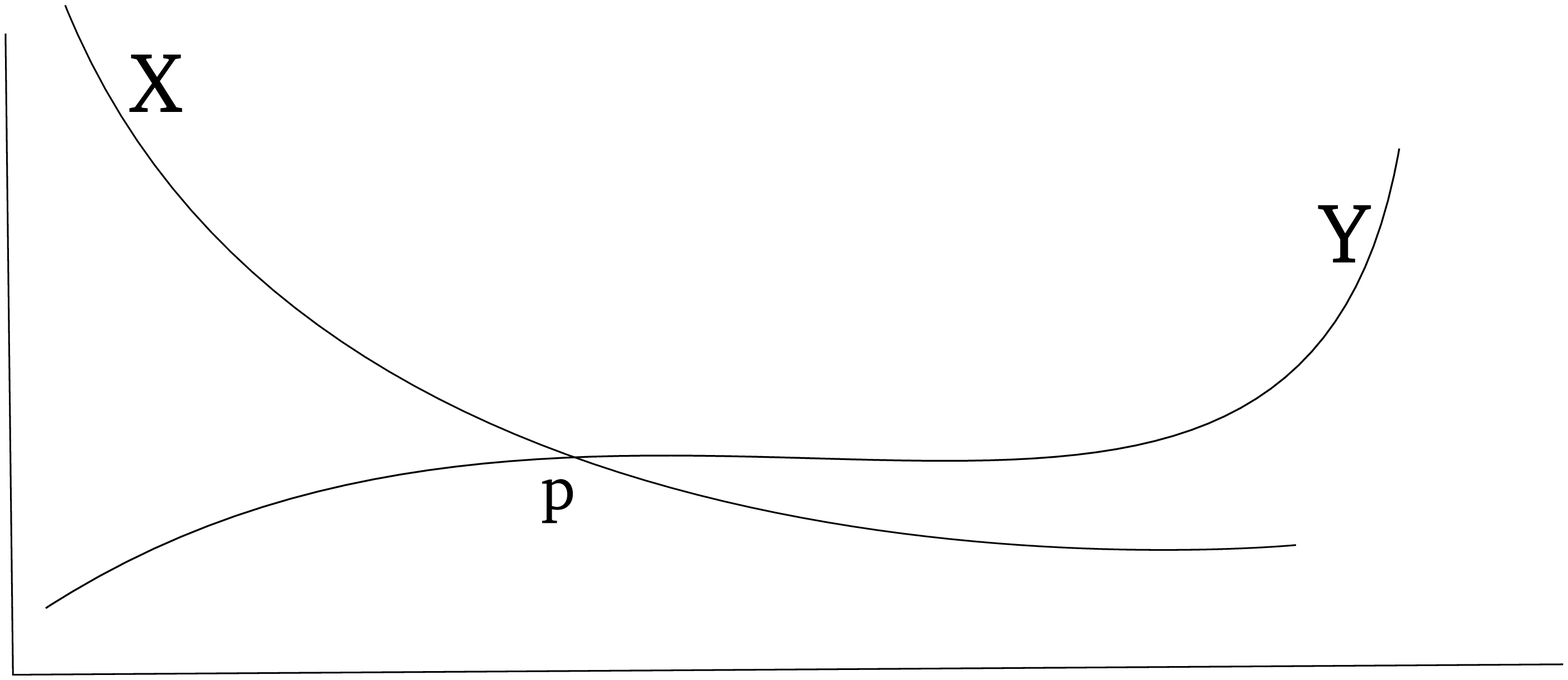}\\
Figure 2:\\$X\prec^+_p Y,Y\prec^-_p X$ \end{center}}}
\end{tabular}
\end{table}
\end{center}
It is easy to see that
``touching'' is a symmetric relation. If $\CM$ expands a real closed
field and $X$ and $Y$ are $C^1$-curves then this notion of tangency
is stronger than the two curves having the same tangent space  at
$p$.

%
%

\begin{defnn} Let $X$ and $Y$ be two definable plane curves.
Assume that $p=\la x_0,y_0\ra\in X\cap Y$ and that $X$ and $Y$ are
graphs of functions, $f(x)$ and $g(x)$, respectively,  near $p_0$.
We say that $X\fleq^+_p Y$ ($X\prec^+_p Y$) if there is $x_1>x_0$
such that $f(x)\leq g(x)$
($f(x)<g(x)$) for all $x$ in the interval $(x_0,x_1)$.

We define $X\fleq^-_p Y$ ($X\prec^-_p Y$) if there exists $x_1<x_0$
such that for all $x\in (x_1,x_0)$, we have $f(x)\leq g(x)$
($f(x)<g(x)$). We write $X\fleq_p Y$ ($X\prec_p Y$) if
$Y\fleq^-_p X$ ($Y\prec^-_p X$) and $X\fleq^+_p Y$ ($X\prec^+_p Y$).
\end{defnn}

If $X$ and $Y$ are graphs of $C^1$ functions $f(x)$, $g(x)$, with
respect to a definable real closed field, and $f(x_0)=g(x_0)$,
$f'(x_0)<g'(x_0)$ then we have $X\prec_p Y$ (see Figure 2)

 Notice that, in the
setting above, $X$ and $Y$ touch each other at $p$ if and only one
of the following holds:
\\(i) $X\fleq^{-}_p Y$ and  $X\fleq^+_p Y$ (namely, $Y$ touches $X$
 from ``below'') , or
\\(ii) $Y\fleq^-_p X$ and $Y\fleq^+_p X$ ($Y$ touches $X$ from
``above'').

\subsection{Normal families of curves}

\begin{defnn}{\em A definable family of curves $\CF=\{C_q:q\in Q\}$ in $ M^n$
 is called  {\em normal} if for every
  $q_1, q_2\in Q$, the curves $C_{q_1}$ and $C_{q_2}$ intersect
 in at most finitely many points. We say that {\em a normal
 $\CF$ has dimension $k$} if $\dim Q=k$.

A family $\CF$ as above is called {\em almost normal} if for every
$q_1\in Q$ there are at most finitely many $q_2\in Q$ such that
$\C_{q_1}\cap \C_{q_2}$ is infinite. The dimension of $\CF$ is by
definition $\dim(Q)$.
 }
\end{defnn}

 Note that if $\CF$ is almost normal and $q\in Q$ is
generic, there exists a neighbourhood $U\ni q$ such that $C_q\cap
C_{q'}$ is finite for all $q'\in U$. If we take $U$ definable over
generic parameters, then we get by the genericity of $q$, an open
$V\ni q$ such that for all $q_1,q_2\in V$, $C_{q_1}\cap C_{q_2}$ is
finite. Therefore, $\CF_V:=\{C_q:q\in V\}$ is normal. Since the
results of this section are all local in nature, using the above
observation, they are all true for almost normal families (with the
same proofs). For the sake of clarity, the statements are only given
for normal families.

\begin{lemma}\vlabel{basicprop} Let $\CF=\{C_q:q\in Q\}$ be a definable,
$k$-dimensional normal family of curves in $P\sub M^n$, with $\dim
P=2$ and $k>1$. Assume that $\CF$ is definable over $A$.
\\(i) For every $\la p,q\ra\in P\times Q$, if $p\in C_q$ and $\dim(p,q/A)=k+1$ then
$p$ and $q$ are generic over $A$ in $P$ and $Q$, respectively, and
$\dim(p/qA)=1$.
\\(ii) If $P\sub M^2$, $p=\la x_0,y_0\ra$ and $\Mdim(p,q)=k+1$ then there
exist open intervals $I$ and $J$, around $x_0$ and $y_0$, respectively,
such that $C_q\cap I\times J$ is the graph of a continuous, strictly
monotone function.
\end{lemma}
\pf (i) Because for every $q\in Q$ we have $\dim(C_q)=1$, it follows
that  $\dim(q/A)=k$, and hence $\dim(p/qA)=1$. Because the
intersection of any two $C_q$ is finite we have $\dim(p/A)=2$.
\\(ii) Because $p$ is generic in $C_q$, the curve $C_q$ is the
graph of a continuous function $f_q$ near $p$ (either in $x$ or in
$y$). If this function were locally constant, say as a function of
$x$, then, by the genericity of $q$ we would get a $k$-dimensional set of $f_q$
locally constant near $x_0$. Because $k>1$ we would get an infinite set of $q'$
 such that $f_{q'}$ all agree on a whole interval,
contradicting normality.\qed

When $\CF$ is a fixed family of curves parameterised by $Q$ then we
sometimes write $q_1\fleq_pq_2$ instead of $C_{q_1}\fleq_p C_{q_2}$.

\subsection{Nice families of curves}

\begin{lemma}\vlabel{very normal} Let
$\CF=\{C_q:q\in Q\}$ be a $\0$-definable two-dimensional normal family
of curves in $P\sub M^2$. Assume that $p_0\in \Cq$ and that
$\Mdim(p_0,q_0/\emptyset)=3$.
 Then there exist a neighbourhood $U\times W\sub M^2\times Q$
of $\la p_0,q_0\ra$ such that for all $q_1,q_2\in W$ the curves
$C_{q_1}$ and $C_{q_2}$ intersect at most once in $U$ and they do
not touch each other at their point of intersection.
\end{lemma}
\proof Assume that $p_0$ and $q_0$ are as above. Then, by possibly
shrinking $P$ and $Q$ we may assume that for every $q\in Q$ the
curve $C_q\sub P$ is the graph of a continuous partial function
$f_q$ from $M$ into $M$ and that $f_q$ vary continuously with $q$.
\\

 \noindent  {\bf Claim 1} For all $q\neq q_0$ sufficiently close
to $q_0$ the curve $C_q$ does not touch $C_{q_0}$ at $p_0$.
\\

 Indeed, if the claim fails
then there is a 1-dimensional set $Q'\sub Q$ containing $q_0$ and
definable over generic parameters, such that $C_{q'}$ touches $C_{q_0}$ at $p_0$ for every $q'\in Q'$. Because $q_0$ is generic, by standard dimension considerations, we can obtain a definable, one-dimensional set $Q_{p_0}$ such that for
any  $q_1,q_2\in Q_{p_0}$, the curves $C_{q_1}, C_{q_2}$, touch each
other at $p_0$. Now, by varying $p_0$, we may find an open set $U$
containing $p_0$ and a relatively open $Q_1\sub Q$ such that for
every $p\in U$ and every $q_1,q_2\in Q_1$, if $p\in C_{q_1}\cap
C_{q_2}$ then $C_{q_1}$ and $C_{q_2}$ touch each other at $p$. We
may assume that $U$ and $Q_1$ are 0-definable (by choosing them to
be definable over sufficiently generic and independent parameters).
Let us show that this contradicts the fact that $Q$ is
two-dimensional and $\CF$ is normal:

Let $p_0=\la x_0,y_0\ra$ and $Q_{0}=\{q\in
Q_1:f_{q}(x_0)=y_0\}$. Choose $x_1>x_0$ generic (over all mentioned
parameters) and $y_1=f_{q_0}(x_1)$ such that $\la x_1,y_1\ra\in
U$. Since $y_1$ is generic in $M$ over $p_0q_0$, there are
$q_1,q_2\in Q_{p_0}$ such that $f_{q_1}(x_1)<y_1<f_{q_2}(x_1)$. If
we now take any $q\in Q_1$ such that $f_q(x_0)\neq y_0$ then
necessarily $f_q(x_1)\neq y_1$, for otherwise the graph of $f_q$
will have to cross (and therefore not touch) either
$C_{q_1}$ or $C_{q_2}$ in the interval $(x_0, x_1)$. By the
normality of $\CF$ and the genericity of
$(x_0,y_0), (x_1,y_1)$ there are at most finitely many $q\in Q_1$
such that both points are in $C_q$. So there are at most finitely
many $q\in Q_1$ such that $\la x_1,y_1\ra \in C_{q_1}$, in
contradiction to the fact that $\CF$ is two-dimensional. End of
Claim 1.
\\

\noindent{\bf Claim 2} There is an open $U\sub P$ and an open $W\sub
Q$ such that for every $q_1\neq q_2\in U$ and for every $p\in
C_{q_1}\cap C_{q_2}\cap W$, the curves $C_{q_1}$ and $C_{q_2}$ do
not touch each other at $p$.
\\

Indeed, by Claim 1 there is a relatively open $W'\sub Q$ such that
for every $q\neq q_0$ in $W'$, if $p_0\in C_q\cap C_{q_0}$ then
$C_q$ and $C_{q_0}$ do not touch each other at $p_0$. By choosing
the parameters defining $W'$ sufficiently independent, we may assume
that $W'$ is 0-definable.  By the genericity of $q_0$ in $Q_{p_0}$,
we may find an open $W\sub Q$,  $q_0\in W\sub W'$,  such that for
every $q',q''\in W$, if $p_0\in C_{q'}\cap C_{q''}$ then the two
curves do not touch each other at $p_0$. Again, we assume that $W$
is 0-definable. Finally, we may use the
genericity of $p_0$ to obtain the desired $U$. End of Claim 2.

The rest of the proof is extracted from Lemma 4.3 in \cite{PeSt2}, so
we only sketch the argument.

We first may assume that $Q\sub M^2$. By Claim 2,
we may assume that for all $p\in U$ and $q_1, q_2\in Q$, if
$C_{q_1}, C_{q_2}$ intersect in $p\in U$ then either $q_1\fleq_p
q_2$ or $q_2\fleq_p q_1$.

 Using Claim 2 and o-minimal-type  dimension arguments, we may find open definable
neighbourhoods $W$ of $q_0$ and $U$ of $p_0$ such that for every
$p\in U$ and for every $q_1=\la a_1,b_1\ra, q_2=\la a_2,b_2\ra\in W$
such that  $p\in C_{q_1}\cap C_{q_2}$, we have: $q_1\fleq^+_p q_2$
if and only if $a_1\leq a_2$. In particular, the
$\fleq^+_p$-ordering depends only on $q_1, q_2$ and not on $p$.

 We claim that these $U$ and $W$ are the desired neighbourhoods.
Indeed, if $C_{q_1}$ and $C_{q_2}$ intersect more than once in $U$,
for $q_1,q_2\in W$, then take two consecutive points of
intersections $p_1$ and $p_2$. Notice that if $q_1\fleq^+_{p_1} q_2$
then, because we assumed that $C_{q_1}$ and $C_{q_2}$ do not touch
each other, we necessarily have $q_2\fleq^+_{p_2} q_1$. But this
contradicts the fact that the $\fleq^+_p$-ordering depends only on
$q_1$ and $q_2$ (and not on $p$).\qed
\\

The following is a variation of a similar definition from
\cite{PeSt2}.
\begin{defnn} Let $\CF=\{C_q:q\in Q\}$ be a definable normal family of plane
curves,
all contained in $U\sub M^2$, where $Q$ is an open subset of $M^2$.
We say that $\CF$ is {\em a nice family} if the following hold:
\\(i) For every $q_1,q_2\in Q$, $C_{q_1}$ and $C_{q_2}$ intersect
at most once in $U$.
\\(ii) For every $q\in Q$, the curve $C_q$ is the graph of a
partial function $f_q:M\to M$.
\\(iii) The partial function $F(a,b,x)$ which sends $\la a,b,x\ra\in Q\times M$
 to $f_{\la a,b\ra}(x)$ is
continuous in all variables and strictly monotone in each of its
variables. Moreover, in each of the variables, $F$ has the same
monotonicity behaviour, as the other two variables vary. E.g., if
$F(a_1,b_1,-)$ is strictly increasing in the  last variable, for
some $\la a_1,b_1\ra\in W$ then for all $\la a,b\ra\in W$, the
function $F(a,b,-)$ is strictly increasing.
\end{defnn}

\begin{cor}\vlabel{nice} Let $\CF=\{C_q:q\in Q\}$ be a 0-definable normal family of
plane curves, where $Q$ is a two-dimensional subset of $M^2$. Assume
that $p_0\in \Cq$, with $\Mdim(p_0,q_0/\emptyset)=3$.  Then there is
an open neighbourhood $ U\times W$ of $\la p_0,q_0\ra$ such that the
family $\CF_{U,W}=\{C_q\cap U:q\in W\}$ is nice.
\end{cor}
\pf Clause (i) of the definition of a nice family follows from
Lemma \ref{very normal}, while the other clauses follow from
genericity.\qed

\begin{defnn}{\em Given a nice family of plane curves $\CF=\{C_q:q\in Q\}$ and a
definable curve $X\sub M^2$ we say that {\em $X$ is $\CF$-bounded at
$p$} if  there are $q_1, q_2\in Q$ such that $p\in C_{q_1}, C_{q_2}$
and $C_{q_1}\fleq_p X\fleq_p C_{q_2}$ (see figure)}
\end{defnn}

\begin{center}
\includegraphics[height=1.2in,width=2in]{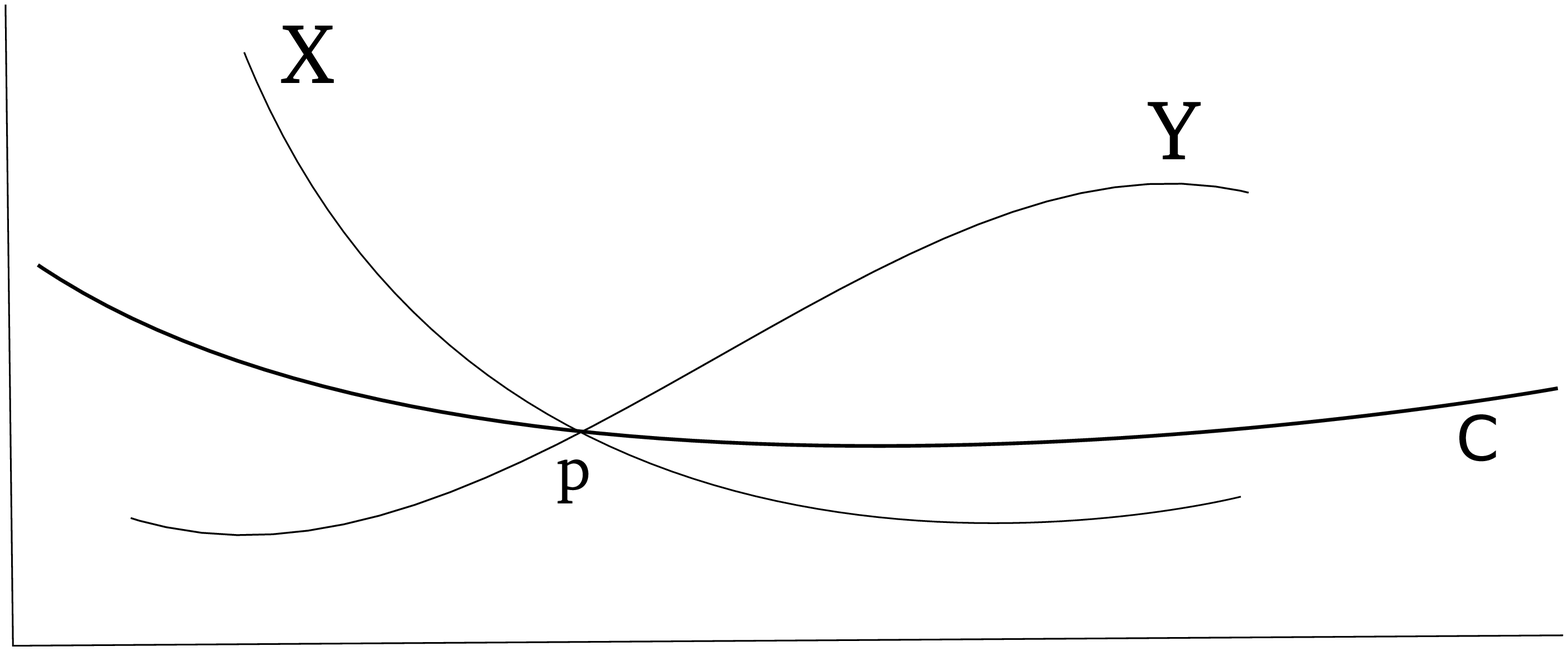}

$C$ is bounded at $p$ by $X$ and $Y$ 
\end{center}

 The following
is just Theorem 10.5 from \cite{PeSt2}. The proof there uses no more
than our definition of a nice family.

\begin{fact}\vlabel{10.5} Let $\CF$ be a nice family of curves in $U\sub M^2$
and let $X\sub U$ be a definable curve that is $\CF$-bounded at
every point. If $p$ is generic in $X$ then there exists a unique
curve from $\CF$ which touches $X$ at $p$.

In particular, if $X=C_q$ then the only curve from $\CF$ that
touches $C_q$ at $p$ is $C_q$ itself.
\end{fact}

\begin{cor}\vlabel{cor1} Let $\CF$ be a two-dimensional normal family of plane curves.
 Let
$X\sub M^2$ be an  $A$-definable curve, $p$ generic in $X$ (over
$A$) and also generic in $M^2$ (over $\emptyset$). Then there are at
most finitely many curves from $\CF$ which touch $X$ at $p$.
\end{cor}
\pf  If there were infinitely many touching curves in $\CF$ we
could find $q\in Q$ such that $\dim(p,q/A)=3$ and $C_q$ touches $X$ at
$p$. By \ref{basicprop}, $q$ is generic in $Q$ over $A$. By Lemma
\ref{nice}, there is a neighbourhood $U\times W$ of $\la p,q\ra$ such
that $\CF_{U,W}=\{\la p, q\ra\in U\times W:p\in C_q\}$ is nice.
Furthermore, we can choose $U$ to be definable over independent
parameters. It follows from \ref{10.5} that $q\in \Mdcl(Ap)$,
contradiction.\qed

\begin{defnn}{\em
Let $\CF$ be a normal family of curves and $X$ a definable curve in
$M^2$. For $p$ a $C^0$-point of $X$, let us denote by $\tau(X,p)$
the set of all $q\in Q$ such that $C_q$ touches $X$ at $p$.

 We let
$$\tau(X)=\bigcup \{\tau(X,p):p \mbox{ a $C^0$-point of } X\}\sub Q.$$

We say that $X$ is {\em $\CF$-linear near $p$} if there exist a
neighbourhood $U$ of $p$ and $q\in Q$ such that for every $p'\in
U\cap X$, we have $\tau(X,p')=\{q\}$.}
\end{defnn}

If $\CF$ is a nice family then $X$ is $\CF$-linear near a generic
$p$ if and only if $X$ equals near $p$ to some $C_q$ in $\CF$. If
$\CF$ is not assumed to be nice then we just have the ``only if''
direction. Namely, it is possible that $C_q\in \CF$ has curves from
$\CF$ (other than itself) touching it at $p$, in which case $X=C_q$
will not be $\CF$-linear.

%

Assume that $\CF$ is a nice
family of curves,  $X$ is $\CF$-bounded and not $\CF$-linear near
any of its points. Then for any $p\in X$ for which $\tau(X,p)\neq
\emptyset$,  $\tau(X,p)$ and $p$ are inter-definable in $\CM$ over
$A$. In particular, $\tau(X)$ is a definable 1-dimensional set.
 By corollary \ref{cor1} this would still be true for normal $\CF$
  provided that there is no  $p\in X$ such that $\{q\in Q: p\in C_q\}$
   is two-dimensional. Furthermore, by \ref{cor1}, if $p\in X$ is
    generic and $q\in \tau(X,p)$ then $q$ is generic in $\tau(X)$ over $A$.

\bigskip

\noindent{\bf Infinitesimals } It is useful for much of what is
coming next to introduce the notion of infinitesimals: Given $p\in
M^n$, we denote by $\nu_p$ the {\em infinitesimal neighbourhood of
$p$}, defined as follows: Given $\CM^*$ an $|M|^+$-saturated
elementary extension  of $\CM$, $\nu_p$ is the intersection of all
$\CM$-definable open neighbourhoods of $p$ in $\CM^*$. Notice that
for $a=\la a_1,\ldots, a_n\ra \in M^n$, we have
$\nu_a=\nu_{a_1}\times\cdots\times\nu_{a_n}$.

Even though the set $\nu_p$ is not definable, most of the statements
involving $\nu_p$ appearing in this paper can be restated in a first order
manner. For example, the statement ``For every $q\in \nu_{q_0}$ the
set $X\cap C_q\cap \nu_{p_0}$ is nonempty'' is equivalent to ``There exist
  neighbourhoods $U$ of $p_0$ and $W$ of $q_0$ such that for all $q\in W$,
   the set $X\cap C_q\cap U$ is nonempty''.

\begin{lemma}\vlabel{tangency2}
Let  $\CF$ be a nice family of plane curves in $U\sub M^2$, $X\sub
U$ an $A$-definable curve that is $\CF$-bounded. Assume that $p_0\in X$
is generic over $A$, that $X$ is not $\CF$-linear
near $p_0$ and that $\Cq$ touches $X$ at $p_0$. Let $W_1, W_2$ be
the locally definably connected components of $M^2\setminus \tau(X)$
at $q_0$.

Then, there is $i\in \{1,2\}$, such that the following are
true:
\begin{enumerate}
\item
Given any $q\in W_i\cap \nu_{q_0}$, the set $C_{q}\cap X\cap
\nu_{p_0}$ is empty.
\item Given any $q\in W_{3-i}\cap \nu_{q_0}$,
the set $C_{q}\cap X\cap \nu_{p_0}$ contains at least two points.
\end{enumerate}
\end{lemma}

\proof Write $p_0=\la x_0,y_0\ra$, $q_0=\la a_0,b_0\ra$. As discussed above,
 $q_0\in\tau(X)$ is generic over $A$ so in
particular it is a $C^0$ point of $\tau(X)$. Consider the continuous
function $F(a,b,x)=f_{\la a,b\ra}(x)$, as given by the definition of
nice families. Without loss of generality the domain of $F$ is a
product of intervals $I_1\times I_2\times I_3$ and its range is
contained in the interval $I_4$.
 Since $p_0$ is generic
in $X$, the set $X$ itself is also locally, near $p_0$, the graph of a function,
and without loss of generality we assume it is a function of $x$.
Call this function  $g(x)$. (If $X$ is not a function of $x$, then,
by genericity it is locally of the form $x=c$; in this case, notice
that since  $\CF$ is nice, we can interchange the roles of $x$ and
$y$ in $\CF$ and consider it as a family of functions of $y$).

 There are several cases to consider, but we will handle only one of them
(the rest can be handled in a similar way). We assume that $\Cq$
touches $X$ from above at $p_0$.
 Namely, for all $x$ near $x_0$, we have $f_{q_0}(x)\geq g(x)$.
Because $p_0$ is generic in $X$ and $\Cq$ is the unique touching
curve from $\CF$ at $p_0$, there is a neighbourhood $U$ of $p_0$ such
that for all $p\in U\cap X$, there is a touching curve $C_q\in \CF$
to $X$ at $p$ and $C_q$ touches $X$ from above there.

Since $q_0$ is generic in $\tau(X)$ over $A$, we may assume that
$\tau(X)$ is the graph of a function in the first variable.

Let us assume (by the definition of nice families) that $F(a,-,x)$
is strictly increasing in the second variable, near $b_0$, for all
$a,x$ close to $a_0, x_0$, respectively. Fix an interval $I$
containing $\nu_{x_0}$ such that for all $\la a,b\ra \in \nu_{q_0}$
we have $I\sub dom(f_{a,b})$. It follows that for all $b>b_0$ and
$x\in I$, we have $F(a_0,b,x)>F(a_0,b_0,x)>g(x)$. In particular, the
curve $C_{a_0,b}$ lies above the curve $X$ in some neighbourhood of
$p_0$ and therefore the two curves do not intersect near $p_0$ (see
Figure 3, with $C_{q_0}$ for $C_{a_0b}$ ).

\begin{table}[h]
\begin{center}%
\begin{tabular}
{ll}%
{\parbox[b]{2in}{
\begin{center}
\includegraphics[height=1.2in,width=2in]{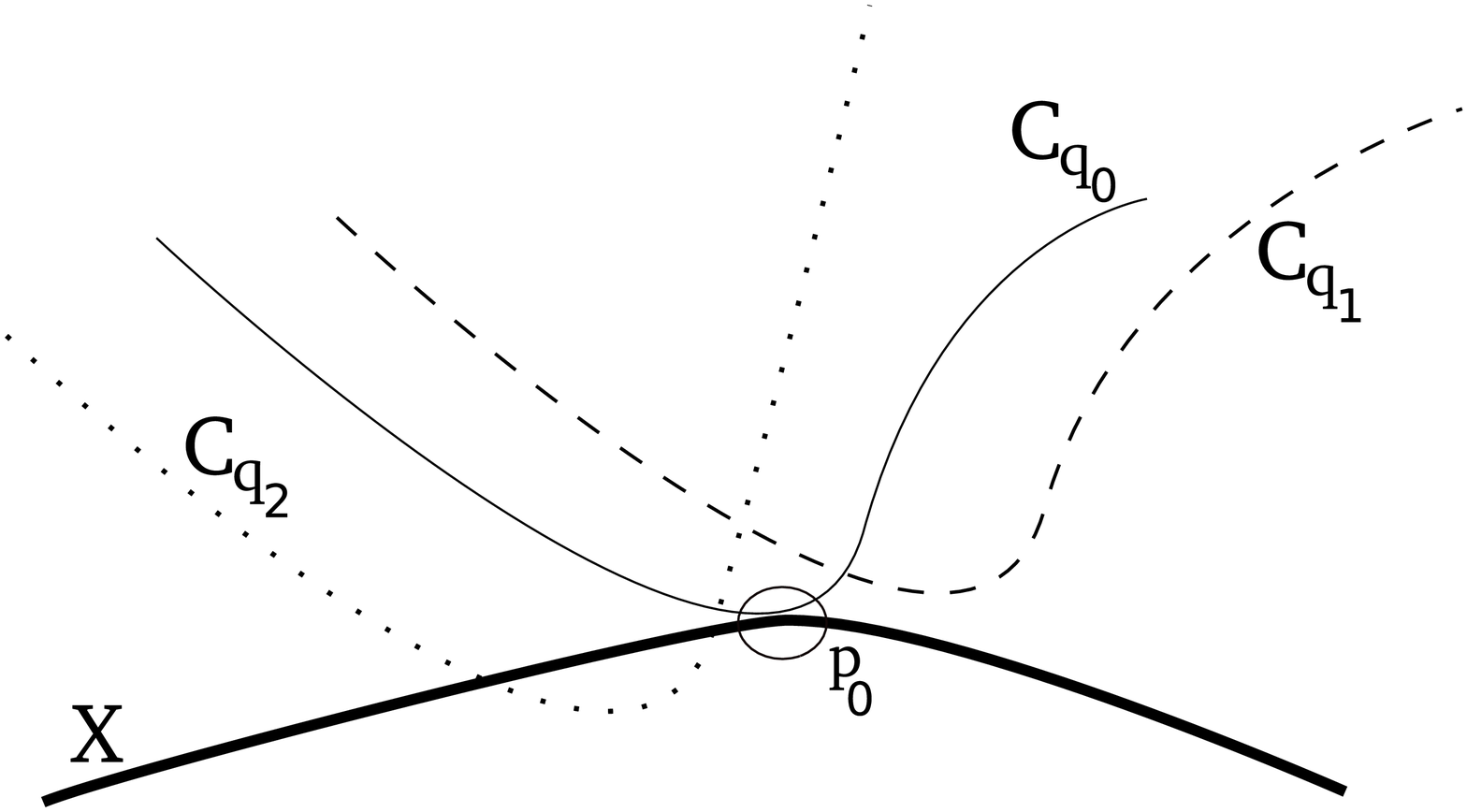}\\
Figure 3:\\ The $P$ universe \end{center}}}
{\parbox[b]{4in}{\begin{center}
\includegraphics[height=1.2in,width=2in]
{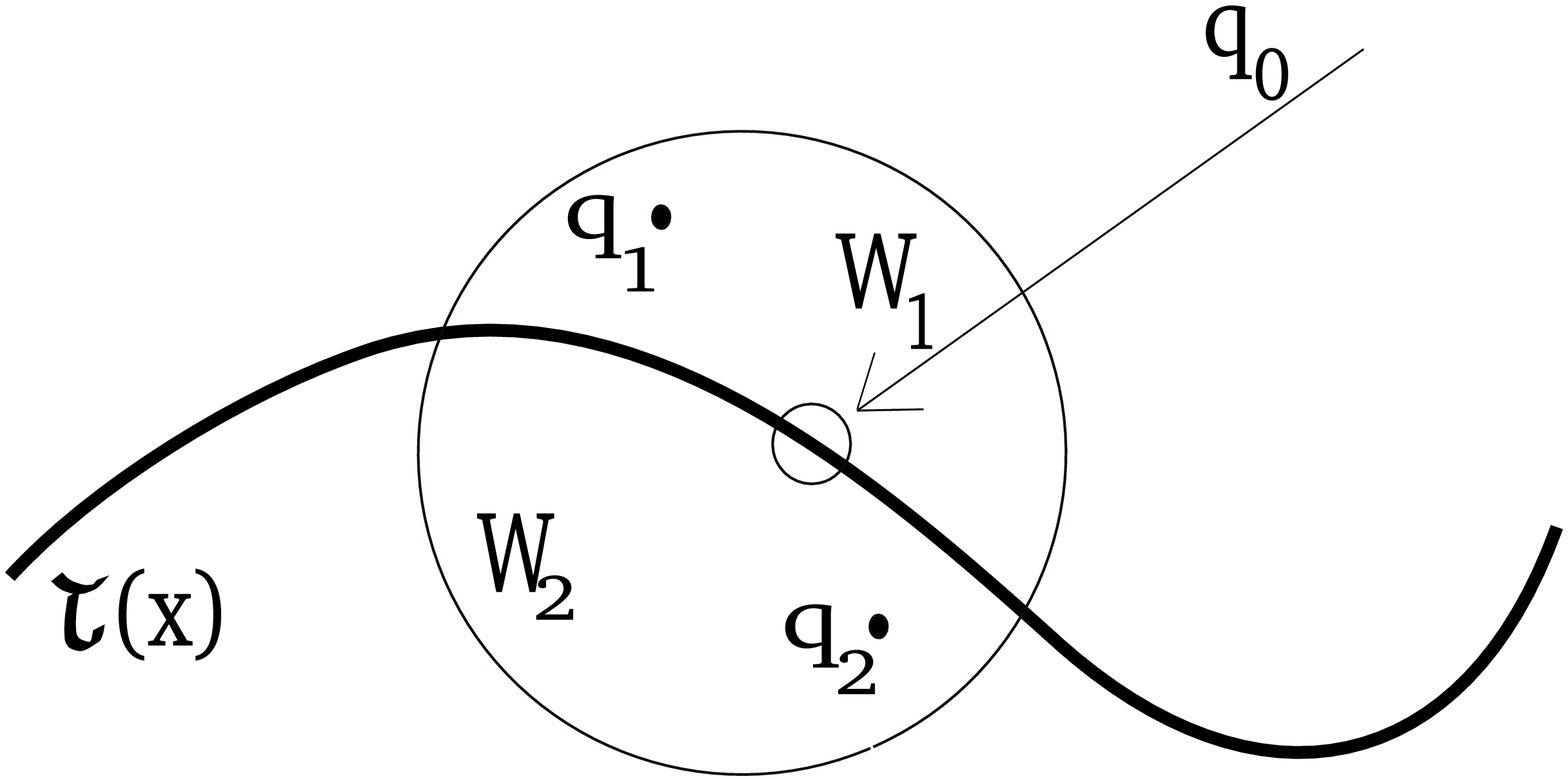}\\
Figure 4:\\The $Q$ universe \end{center}}}%
\end{tabular}
\end{center}
\end{table}

Now, as we pointed out above, for all $q=\la
a,b\ra\in \tau(X)\cap \nu_{q_0}$, the curve $C_q$ touches $X$ from
above at some $p\in \nu_{p_0}$. If $q'=\la a,b_1\ra\in \nu_{q_0}$
and $b_1>b$ then by the same argument $C_{q'}\cap X\cap
\nu_{p_0}=\emptyset$.

 So, in the case we just considered, the local
component $W_1$ of $M^2\setminus \tau(X)$ which lies above the curve
$\tau(X)$ satisfies: For all $q_1\in W_1\cap \nu_{q_0}$,
$C_{q_1}\cap X\cap \nu_{p_0}=\emptyset$  (see Figures 3,4, with
$C_{q_1}$ for $C_{a_0b_1}$).

Let us consider the local component $W_2$ which lies below
$\tau(X)$. If we take $b_1<b_0$, and $\la a_0,b_1\ra\in \nu_{q_0}$
then
\begin{equation}\label{eq1}
f_{a_0,b_1}(x_0)<f_{a_0,b_0}(x_0)=g(x_0).\end{equation}

At the same time, by the previous argument there are $x_1<x_0<x_2$
both in $\nu_{x_0}$ such that $f_{a_0,b_0}(x_1)>g(x_1)$ and
$f_{a_0,b_0}(x_2)>g(x_2)$. Therefore, by continuity considerations, if $b_1$ is
sufficiently close to $b_0$ then there are
$x_1',x_2'\in \nu_{x_0}$, $x_1'<x_0<x_2'$ such
that
\begin{equation}\label{eq2} f_{a_0,b_1}(x_1')>g(x_1')\,\,\mbox{ and }\,\,
f_{a_0,b_1}(x_2')>g(x_2').\end{equation}

It follows by continuity from (\ref{eq1}) and (\ref{eq2}) that for
some $x_1'', x_2''\in \nu_{x_0}$ we have
\begin{equation}\label{eq3} f_{a_0,b_1}(x_1'')=g(x_1'')\,\,\mbox{ and }
f_{a_0,b_1}(x_2'')=g(x_2'').\end{equation} Hence, $X$ and
$C_{a_0,b_1}$ intersect at least twice in $\nu_{p_0}$, for $b_1<b$
and close to $b_0$ (see, again, Figures 3 and 4).

Repeating the same argument for every $q=\la a,b\ra\in
\tau(X)\cap \nu_{q_0}$ and $\la a,b_1\ra\in \nu_{q_0}$ with $b_1<b$,
we conclude that $|C_{a,b_1}\cap X\cap\nu_{p_0}|>1$.\qed

\begin{lemma}\vlabel{bounded} Let $\CF$ be a nice family.
Let $X\sub M^2$ be an $A$-definable curve, $p_0\in X$ generic over
$A$. If there is a curve from $\CF$ that touches $X$ at $p_0$ then
$X$ is $\CF$-bounded near $p_0$.
\end{lemma}
\proof Assume that $\Cq$ touches $X$ at $p_0$. By Lemma
\ref{tangency2} and its proof, there are $q'$ arbitrarily close to $q_0$
and $p_1,p_2\in \nu_{p_0}$ such that $C_{q'}\fleq_{p_1} X$ and
$X\fleq_{p_2} C_{q'}$. But then the two sets $$\{p'\in X:\exists
q'\,\, C_{q'}\fleq_{p'} X\}$$ and
$$\{p'\in X:\exists q'\,\, X\fleq_{p'} C_{q'}\}$$ intersect $\nu_{p_0}$
nontrivially. Since these sets are definable just over the
parameters defining $X$ and $\CF$, and since $p_0$ is generic in $X$
over these parameters, the sets must contain all of
$X\cap \nu_{p_0}$, hence $X$ is
$\CF$-bounded in a neighbourhood of $p_0$.\qed

\begin{lemma}\vlabel{transversal} Let $\CF$ be a nice family of curves, $X$
an $A$-definable curve, $p\in X$ generic over $A$ and $q_0\in Q$.
If $p\in C_{q_0}\cap X$ then one and only one of the following
holds:
\\(i)  $C_{q_0}$ touches $X$ at $p$, or
\\(ii) For every $q\in \nu_{q_0}$, we have $|C_{q}\cap X|=1$.
\end{lemma}
\proof Notice that by \ref{tangency2}, clauses (i) and (ii) cannot
hold simultaneously. Assume that $\Cq$ does not touch $X$ at $p$.

 Then $\Cq$ intersects both local components of
$M^2\setminus X$ in $\nu_p$. By continuity, $C_{q'}$ intersects both components of $M^2\setminus X$ in $\nu_p$ for all $q'\in
\nu_{q_0}$. By the definable connectedness of $C_{q'}$ (or
equivalently, the continuity of the function $f_{q'}$), every such curve
$C_{q'}$ must intersect $X$ in $\nu_p$ at least once.

Assume now, towards contradiction, that there exists $q\in
\nu_{q_0}$ such that $C_{q}\cap X\cap \nu_p$ contains more than one
point. We claim that $X$ is $\CF$-bounded near $p$.  There are two
possibilities to consider:

(i) $C_q$ touches $X$ at one or more of these points of
intersection. In this case, by Lemma \ref{bounded} above, $X$ is
$\CF$-bounded in some open subset of $\nu_p$. Because of the
genericity of $p$ (over the parameters defining $\CF,X$), the curve
$X$ is $\CF$-bounded in a neighbourhood of $p$.

(ii) $C_{q}$ does not touch $X$ at any of their points of
intersection. In this case, as we observed before, we can find
$p_1,p_2\in \nu_p$ such that $C_{q}\fleq_{p_1} X$ and $X\fleq_{p_2}
C_{q}$. This was already shown to imply that $X$ is $\CF$-bounded in
$\nu_p$.

Since $X$ is $\CF$-bounded, by \ref{10.5} there is a unique $q'\in
Q$ such that $C_{q'}$ touches $X$ at $p$. However, if we now replace
$Q$ by a smaller open neighbourhood of $q_0$ then the above argument
shows that $X$ is still $\CF$-bounded near $p$ with respect to the
restricted  $\CF$ and therefore there is still a curve in this
family that touches $X$ at $p$. This curve must therefore be $\Cq$,
contradicting the assumption that $\Cq$ does not touch $X$ at $p$.
\qed

The last lemma can be seen as saying that in a nice family of
curves, every  curve that intersects $X$  at a generic point either
touches $X$ or is intersects it transversally.

\begin{lemma}\vlabel{specialcase1}
Let $\CF_0=\{C_q:q\in Q_0\}$, $\CF_1=\{D_q:q\in Q_1\}$  be two
0-definable two-dimensional nice families of plane curves. Let
$q_0,q_1$ be generic in $Q_0$, $Q_1$, respectively, and $p_0\in
C_{q_0}\cap D_{q_1}$ such that
$\Mdim(p_0,q_0/\emptyset)=\Mdim(p_0,q_1/\emptyset)=3$. Assume that
there exists some A-definable curve $X$ such that $p_0$ is also
generic in $X$ over $A$. If $C_{q_0}$ and $D_{q_1}$ both touch $X$
at $p_0$ then $\Cq$ and $D_{q_1}$ touch each other at $p_0$, and $q_0$ and $q_1$
are inter-definable over
$p_0$.
\end{lemma}
\pf Notice that if $C_{q_0}$ and $D_{q_1}$ touch $X$ on opposite
sides, then they clearly touch each other as well (see Figure 5). Assume then that
the two curves touch $X$  "from above". Because of the
 genericity of $p_0$ in $X$, for every $p\in \nu_{p_0}$,
 there is $q\in \nu_{q_0}$ such that $C_q$ touches $X$ from above at $p$.
 Similarly to the proof of \ref{bounded},
  it follows that $D_{q_1}$ is $\CF_0$-bounded near $p_0$ (see Figure 6).
  Moreover, even if we shrink $Q_0$ to a smaller neighbourhood of
  $q_0$, the curve $D_{q_1}$ is still $\CF_0$-bounded with respect
  to this new family. As in the proof of Lemma \ref{transversal}, it
  follows (using \ref{10.5}), that $C_{q_0}$ touches $D_{q_1}$ at
  $p_0$. ($D_{q_1}$ plays the role of $X$ from \ref{10.5}).
\begin{table}[h]
\begin{center}%
\begin{tabular}
{ll}%
{\parbox[b]{2in}{\begin{center}
\includegraphics[height=1.2in,width=2in]{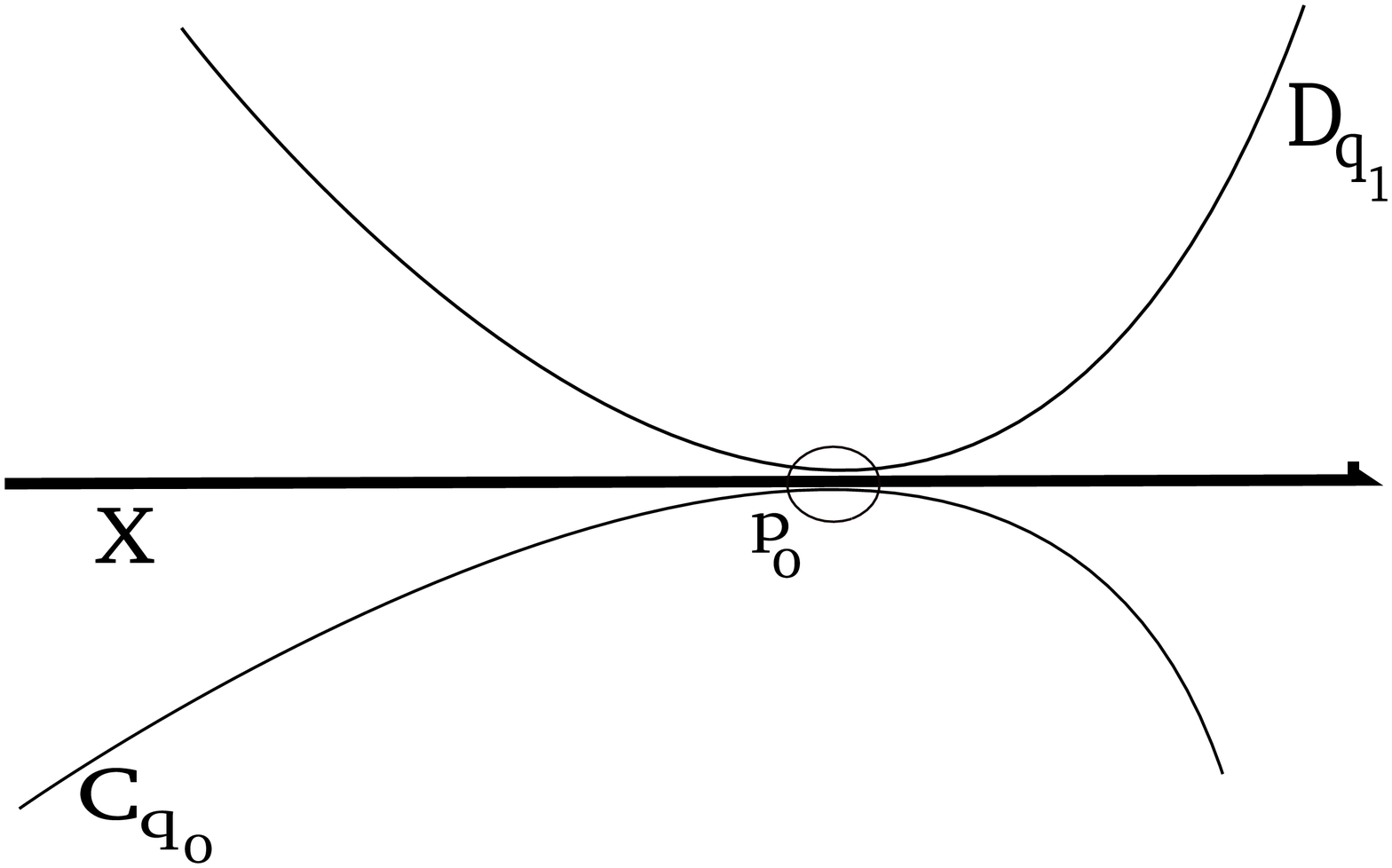}\\
Figure 5: \\$C_{q_0}, D_{q_1}$ on opposite sides \end{center}}}
{
\parbox[b]{4in}{
\begin{center}
\includegraphics[height=1.2in,width=2in]{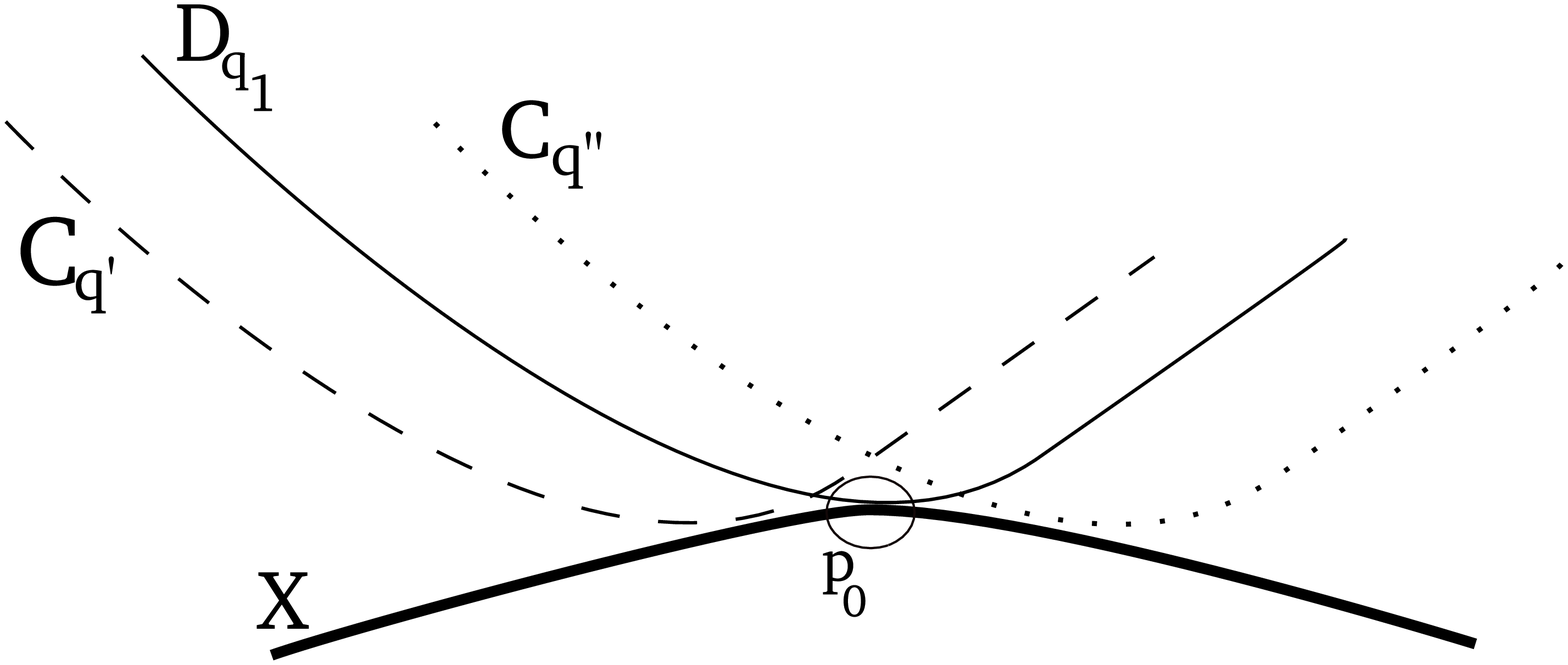}\\
Figure 6: \\ $D_{q_1}$ is $\CF_0$-bounded \end{center}}}
\end{tabular}
\end{center}
\end{table}
   By
 \ref{cor1} (applied to
 the family $\CF_0$ and the curve $D_{q_1}$), if the two curves
 touch each other at $p_0$ then we have
 $q_0\in \Mdcl(q_1,p_0)$. Similarly, $q_1\in dcl(q_0,p_0)$.
 \qed


\subsection{Normal families and differentiability}

{\em Throughout this subsection we assume that $\CM$ expands a real closed field.}

 \begin{lemma}\vlabel{diff1} Let $\{X_t:t\in T\}$ be a 0-definable
normal family of
 plane curves, of dimension at least $2$. Let $t_0\in T$ be generic, $X=X_{t_0}$,
 and $p_0=\la x_0,y_0\ra\in X$ generic over $t_0$. Then $X$ is
 locally the graph of strictly monotone $C^1$-function, $y=f_{t_0}(x)$.
If $d$ is the derivative of $f_{t_0}$ at $x_0$ then
 $\Mdim(d/p_0)=1$.
\end{lemma}
\pf The first part follows from genericity together with the fact
that $\Mdim T\geq 2$ (otherwise, $X$ could be of the form $x=c$). As
for the second part, if $d\in \Mdcl(p_0)$ then there is a definable
function $d(p)$, defined in a neighbourhood $U_0$ of $p_0$, and there exists
a neighbourhood $T_0$ of $t_0$ in $T$,
such that for every $p\in U_0$ and every $X_t$ containing $p$, with
$t\in T_0$, the derivative of $X_t$ at $p$ equals to $d(p)$.

By the Uniqueness Theorem for differential equations in o-minimal
structures (see Theorem 2.3 in \cite{OtPePi}), for every $X_t$,  $t\in
T_0$, if $p_0\in X_t$ then $X_t$ and $X_{t_0}$ coincide in some
neighbourhood of $p_0$. This contradicts normality.\qed

\begin{lemma}\vlabel{diff2} Assume that $\CF=\{C_q:q\in Q\}$
 is a $\0$-definable two-dimensional
normal family of plane curves. If $p_0=\la x_0,y_0\ra$ is generic in
$M^2$ then for all but finitely many $q\in Q$,  if $p_0\in C_q$,
then  $C_q$ is the graph of a $C^1$-function $y=f_q(x)$ near $p_0$.
For every such $q$, $f'_q(x_0)$ and $q$ are inter-definable over
$p_0$.
\end{lemma}
\pf For $p\in M^2$ let $\ell_{p}=\{q\in Q: p\in
C_q\}$. The first part of the lemma follows from the
fact that if $q$ is generic in $\ell_{p_0}$ over $p_0$ then $p_0$ is
generic in $C_q$ over $q$ (by \ref{basicprop}).

Assume now $f_q(x)$ is $C^1$ near $x_0$ and that $d=f'_q(x_0)$. If
there are infinitely many $q'\in \ell_{p_0}$ such that
$f'_{q'}(x_0)=d$ then first of all there is $q_0\in \ell_{p_0}$ such
that $\Mdim(q_0/p_0)=1$ and $d=f'_{q_0}(x_0)$. Furthermore,  because
$\Mdim (\ell_{p_0})=1$, we have $d\in \Mdcl(p_0)$. This contradicts
Lemma \ref{diff1} with $C_{q_0}$ here taken as $X$.\qed

\begin{defnn}{\em  Let $X$ be a definable curve. We say that $X$ has rank
$k$ over $A$ if there is a an A-definable normal family of curves
$\{X_t:t\in T\}$ and $t_0\in T$ generic over $A$, such that
$X=X_{t_0}$ and $\Mdim(T)=k$. We say that {\em $X$ has rank $k$}  if
$A=\emptyset$.}
\end{defnn}

The following is easy to verify:
\begin{rmk}\vlabel{localrank} The notion of rank is well-defined,
even locally. Namely, let $X\sub M^k$ be a definable curve of rank
$k$ over $A$, and $U$ an $A$-definable open set with $U\cap X\neq
\0$. Then for every A-definable normal family of curves $\{X_t:t\in
T'\}$, if there exists $t$ in $T'$ such that $X\cap U=X_t\cap U$
then $\Mdim(T')\geq k$. In particular, $X\cap U$ has also rank $k$
over $A$.
\end{rmk}

We can now establish the connection between the differential notion
of tangency and that of touching curves.
\begin{theorem}\vlabel{diff3} Let $\CF=\{C_q:q\in Q\}$
be an $\0$-definable normal family of plane curves of dimension
greater than $1$. Let $X$ be an $A$-definable plane curve of rank
greater than $1$, $p_0=\la x_0,y_0\ra\in X$ generic over $A$, and
$d$ the derivative of the function associated to $X$ at $x_0$.

If $p_0\in \Cq$ for $q_0\in Q$ and if $f'_{q_0}(x_0)=d$ (in the
notation above) then $\Cq$ and $X$ touch each other at $p_0$ and
$\Mdim(p_0,q_0/\emptyset)=3$.
\end{theorem}
\pf We assume that all sets are 0-definable. By \ref{basicprop},
$p_0$ is generic in $M^2$. It follows from \ref{diff1} that
$\Mdim(d/p_0)=1$, and from \ref{diff2} it follows that
$\Mdim(q_0/p_0)=1$, hence $\Mdim(p_0,q_0/\emptyset)=3$.

We may assume now that $X$ is the graph of a $C^1$-function $y=h(x)$
near $x_0$. By Theorem \ref{nice}, we may assume that $\CF$ is a
nice family.

Assume towards contradiction that $\Cq$ does not touch $X$ at $p_0$.
Say, we have $X\prec_{p_0} \Cq$. It follows from continuity of the
family that for every $q\in \nu_{q_0}$, there is $x>x_0$, $x\in
\nu_{x_0}$, such that $f_q(x)> h(x)$. By \ref{diff1}, we can find
$q\in \ell_{p_0}\cap \nu_{p_0}$ such that $f'_q(x_0)<d$ and
therefore $C_q\fleq_{p_0} X$. These two last facts, together with
the continuity of $f_q$, imply that there is $x>x_0$, $x\in
\nu_{x_0}$, such that $f_q(x_0)=h(x_0)$. In particular, $|C_q\cap
X\cap \nu_{p_0}|>1$, contradicting \ref{transversal}.\qed

\subsection{The Duality Theorem }

For $\CF=\{C_q:q\in Q\}$ a nice family of curves in $P\sub M^2$,
with $Q\sub M^2$, we let $\CL=\{\ell_p:p\in P\}$ be the dual family,
defined by $q\in \ell_p \Leftrightarrow p\in \C_q$. It is not hard
to see that the monotonicity and continuity assumptions in the
definition of a nice family imply that for every $p\in P$, either
$\ell_p$ is empty or $\Mdim (\ell_p)=1$.  Furthermore, because
$$q_1,q_2\in \ell_{p_1}\cap \ell_{p_2}
\,\,\Leftrightarrow\,\, p_1,p_2\in C_{q_1}\cap C_{q_2},$$ for
$p_1\neq p_2$ we have $|\ell_{p_1}\cap \ell_{p_2}|\leq 1$.
 Finally,
notice that if a definable function $G(x,y,z)$ is monotone and
continuous in each variable separately then it is continuous as a
function of all variables. We thus have:

\begin{claim} If $\CF$ is a nice family then the dual family $\CL$
is also nice.
\end{claim}

 We can now establish the following  duality:

\begin{theorem}\vlabel{duality} Assume that $\CF$ is a nice
family of curves. Let $X\sub M^2$ be an $\CM$-definable curve and
let $p_0$ be a generic point of $X$ such that $X$ is not
$\CF$-linear near $p_0$. If $\Cq$ touches $X$ at $p_0$ then
$\ell_{p_0}$ touches $\tau(X)$ at $q_0$.
\end{theorem}
Before proving the theorem, here is an example:

 \noindent {\bf Example:} Let $\CF$ be the family of all affine lines in $\reals^2$,
$C_{a,b}=\{y=ax+b\}$. It is easy to see that the dual family $\CL$
is also the family of affine lines. Take $X=\{y=x^2\}$. Then
$\tau(X)=\{\la u,v\ra:v=-u^2/4\}$. Given $p=\la x_0,x_0^2\ra \in X$,
the $\CF$-tangent curve at $p$ is $C_q$, where $q=\la 2x_0,
-x_0^2\ra$. We also have $\ell_p=\{\la u,v\ra: v=-x_0u+x_0^2\}$.

The slope of $\ell_p$ equals everywhere to $-x_0$, which is exactly
the slope of $\tau(X)$ at $q=\la 2x_0,-x_0^2\ra$. Since $q$ clearly belongs
 to $\ell_p$ it follows that
$\ell_p$ is tangent to $\tau(X)$ at $q$.
\\

\begin{proof}[Proof of theorem] Since $p_0$ and $q_0$ are inter-definable
 (over the parameters defining $X$), the point $q_0$ is generic in
 $\tau(X)$ as well.

\noindent{\bf Claim } The curve $\ell_{p_0}$ is not equal to the
curve  $\tau(X)$ in a neighbourhood of $q_0$.

Indeed, notice that if $\ell_{p_0}$ equals to $\tau(X)$ in a
neighbourhood of $q_0$ then the same is true for all $p$ in  a
relatively open $X_1\sub X$ containing $p_0$. But then all of these
$\ell_p$ intersect each other at infinitely many points,
contradicting the normality of $\CL$.

 Let $W_1, W_2$ be the two local components of $M^2\setminus
\tau(X)$ at $q_0$. Without loss of generality, by Lemma
\ref{tangency2}, for every $q\in W_1\cap \nu_{q_0}$, we have
$C_q\cap X\cap\nu_{p_0}=\emptyset$. However, for every $q\in
\ell_{p_0}$ we clearly have $p_0\in \C_q\cap X\cap \nu_{p_0}$.
Hence, near $q_0$, the entire curve $\ell_{p_0}$ is contained in $W_2$, so
 $\ell_{p_0}$ touches $\tau(X)$ at $q_0$.
\end{proof}

Notice that the Duality Theorem implies, in the notation of the
theorem, that if $X$ is $\CF$-bounded near $p_0$ then $\tau(X)$ is
$\CL$-bounded near $q_0$ (see Lemma \ref{bounded}).

\bigskip

\noindent{\bf Remark} After discovering the above duality we found
out that it has a classical analogue in projective geometry. In that
setting a duality exists between an algebraic variety and the
variety of all its tangent hyper-planes. See, e.g., \cite{BrKn}.
\\

We need the following corollary;
\begin{cor}\vlabel{specialcase} Let $\CF=\{C_q:q\in Q\}$ be a
  0-definable two dimensional normal family of curves in $P$, where
  $\dim P=2$. Assume that $q_0$ is generic in $Q$,
and $X\sub M^2$ is an $A$-definable curve. Assume that $p_0, p_1$
are generic elements of $X$ and that $X$ is not $\CF$-linear near
any of them. If  $\Mdim(p_0,q_0/\emptyset)=3$ and $\Cq$ touches $X$
at both $p_0$ and $p_1$ then $p_1\in \Mdcl(p_0,q_0)$ (although $X$
is $A$-definable!).
\end{cor}
\pf We may clearly assume that $p_1\notin \Mdcl(q_0)$. By
\ref{nice}, we may shrink $P$ and $Q$, and obtain two nice families
of curves $\CF_0=\{C^0_q:q\in Q_0\}$  and $\CF_1=\{C^1_q:q\in
Q_1\}$, such that $C^0_{q_0}$ and $C_{q_0}$ agree in a neighbourhood
of $p_0$ and similarly, $C^1_{q_1}$ and $C_{q_1}$ agree in a
neighbourhood of $p_0$. By choosing $\CF_0$ and $\CF_1$ properly we
may absorb all new parameters and still assume that they are
0-definable.

By \ref{duality} the curves $\ell_{p_0}$ and $\ell_{p_1}$ touch
$\tau(X)$ at $q_0$. We now apply \ref{specialcase1} to the families
$\CL_0$ and $\CL_1$ (the dual families of $\CF_0$ and $\CF_1$), the
points $p_0, p_1$ and $\tau(X)$, and conclude that $p_0\in
\Mdcl(q_0,p_1)$. \qed

\bigskip

\noindent{\bf Remarks } \begin{enumerate}
\item Assume that $\CF$ is a $\0$-definable normal
family of plane curves in $P$ (but $\CF$ need not be a nice
family), such that its dual family $\CL$ is also normal. Assume also that $X\sub P$
 is a definable curve of rank $k>1$, $p$
generic in $X$, $\Cq$ touches $X$ at $p_0$ and
$\dim(p_0,q_0/\emptyset)=3$. By working locally, we can apply
Theorem \ref{duality} and conclude that $\ell_{p_0}$ touches
$\tau(X)$ at $q_0$. However, in this case it is possible that there
are finitely many such $C_q$, all touching $X$ at $p_0$, and hence
$\ell_{p_0}$ will touch $\tau(X)$ at each one of these points $q$.

\item In many of the results in this section we  add the extra
assumption that $\Mdim(p,q/\emptyset)=3$. This assumption is indeed
necessary, as the following example shows:

\noindent Let $C=\{\la x,y\ra\in \reals^2 :y=\sqrt{|x|}\}$ and for
$\la a,b\ra\in \reals^2$ let $C_{a,b}=C+\la a,b\ra$. The family
$\{C_{a,b}:\la a,b\ra \in \reals^2\}$ is normal. If $X\sub \reals^2$
is any $C^1$-curve then for any $p_0=\la a,b\ra\in X$, the curve
$C_{a,b}$ touches $X$ at $p_0$, and  clearly $q_0=\la a,b\ra\in
\Mdcl(p_0)$. In this case, if $C_{ab}$ touches $X$ at another point
$p_1$ then it will not be true in general that $p_1\in
\Mdcl(p_0,q_0)=\Mdcl(p_0)$.
\end{enumerate}

Finally, we have:
\begin{lemma}\vlabel{rank} Assume that $\CF$ is a
$\0$-definable nice family of
plane curves, parameterised by $Q$. Let $X$ be a nowhere $\CF$-linear $t$-definable
 plane curve of rank $k>1$
over $A$. For every $p_0\in X$ generic over $A,t$, if $\Cq$ touches
$X$ at $p_0$, then there exists a neighbourhood $W$ of $q_0$ such
that $\tau(X)\cap W$ is a curve of rank $k$ over $A$.
\end{lemma}
\pf Because $q_0$ is generic in $Y=\tau(X)$ and $\CL$ is a nice
family, there are neighbourhoods $U$ of $p_0$ and $W$ of $q_0$ such
that $p\mapsto \tau(X,p)$ induces a bijection of $X\cap U $ and
$\tau(X)\cap W$, and moreover, by the Duality Theorem $\tau(\tau(X\cap U))=X\cap U$.

 Recall that $X=X_t$ for the $A$-definable family $\{X_t:t\in T\}$, where
$\Mdim(t/A)=k$, and consider the family $\{\tau(X_t\cap U) :t\in
T\}$. We may assume that for all $t\in T$ we have $\tau(\tau(X_t\cap
U))=X_t\cap U$.

We claim that this last family must be normal. Indeed, if not then
there are $t_1,t_2\in T$ and an open set $W_1\sub W$ such that
$\tau(X_{t_1}\cap U)\cap W_1=\tau(X_{t_2}\cap U)\cap W_1$. If we now
apply $\tau$ to those sets then, by the Duality Theorem  we will
conclude, for some $U_1\sub U$,  that $X_{t_1}\cap U_1=X_{t_2}\cap
U_1$, contradiction.\qed

\subsection{Limit sets and a theorem of v.d. Dries}

Somewhat surprisingly, the notion of limit sets plays a crucial role
in our  subsequent counting of intersection points of curves. Several
different such notions of  were studied extensively by v.d. Dries 
in \cite{vdDr1}. The definition below, which 
is suitable for an arbitrary o-minimal structure (not necessarily over
the reals) resembles most that of the Hausdorff limit, but is
not restrictred to families of compact, or even bounded sets.

\begin{defnn} {\em Consider a definable set $S\sub M^k\times M^n$.
 For $x\in Z_1=\pi_1(S)$,
 the projection of $S$ on the first $k$-coordinates, let $S_x=\{y\in
 M^n:\la x,y\ra\in S\}.$
Let $\CS=\{S_z:z\in Z_1\}$ be the corresponding family of
subsets of $M^n$. Given $a<b$ in $ M\cup\{\pm \infty\}$, and given a
definable map $\gamma:(a,b)\to Z_1$, we let
$\CS(\gamma)$ be the set of all $y\in M^n$ such that for every open
neighbourhood $V$ of $y$ and
for every $\epsilon \in (a,b)$, there is $t\in (a,\epsilon)$ such
that $(S_{\gamma(t)}\cap V)\neq \emptyset$.

$\CS(\gamma)$ is called  {\em the limit  of $\CS$ along $\gamma(t)$,
as $t$ tends to $a$}. }
\end{defnn}

Notice that equivalent definitions to the above are given by: $y\in
\CS(\gamma)$ iff for every $t>a$ in $\nu_a$, $S_{\gamma(t)}\cap
\nu_y\neq \emptyset$ iff there exists $t>a$ in $\nu_a$ such that
$S_{\gamma(t)}\cap \nu_y\neq \emptyset$.

\noindent Here are some easy facts about limit sets:

(1) If every $S_x$ has dimension $k$ then the dimension of every
limit set as above is at most $k$.

(2) If $\lim_{t \to a} \gamma(t)=z_0$ is a generic point of $Z_1$
(over the parameters defining $\CS$)  then
$\CS(\gamma)=\cl(S_{z_0})=\cl(S)\cap (\{z_0\}\times M^n)$. In
particular, for any other definable $\gamma':(a',b')\to Z_1$ such
that $\lim_{t\to a'}\gamma'(t)=z_0$, we have
$\CS(\gamma)=\CS(\gamma')$.

(3) Assume that $X\sub Z_1$ is a definable one-dimensional set,
$p\in \cl(X)$. Then there are finitely many ways to approach $p$ in
$X$ (depending on the number of local components of
$X\setminus\{p\}$). We write $\CS(X,p)$ for the union of the
finitely many possible limits sets $\CS(\gamma)$ as $\gamma(t)$
tends to $p$ in $X$, and call this  set {\em the limit set of $\CS$
along $X$ at $p$}. If $p$ is a generic point in $X$ then $\CS(X,p)$
is just $\cl(S_p)$.

 The following theorem is due to van den Dries
 (the complete - unpublished - proof
  of which appears in \cite{vdDr2}, and is  similar to the proof of
 the analogues  theorem in \S 9 of \cite{vdDr1}):

\begin{theorem}\vlabel{limit sets} Assume that $\CM$ expands a real
closed field. Given a definable family $\CS$ as above, the family
$$\tilde{\CS}=\{\CS(\gamma)| \gamma:(a,b)\to Z_1 \mbox{ definable }\}$$
 is definable in $\CM$
and $\dim(\tilde{\CS})\leq \dim Z_1$. If furthermore all the
sets in $\CS$ are closed then $\dim(\tilde{\CS}\setminus \CS)<
\dim Z_1$ (By $\dim(\tilde{\CS})$ we mean the smallest possible
dimension for a set of parameters for $\tilde{\CS}$).
\end{theorem}
In this paper we will be using only the first part of the
theorem. Since the proof does not appear elsewhere, we present in
Appendix A, with the author's permission, his proof of the theorem.

\subsection{Special points}

{\em In this section we assume that $\CF=\{C_q:q\in Q\}$ is a
0-definable two-dimensional normal family of curves in $P\sub M^k$,
$\dim P=2$ (but $P$ is not necessarily contained in the plane), such
that its dual family $\CL$ is also normal}.

A few words of explanations are in place regarding these assumptions.
\begin{enumerate}
\item If $\CF$ is an (almost) normal two-dimensional family of curves, its dual family
 need not consist only of curves, i.e. for $p\in P$ the set $\ell_p$ of $\CF$-curves
 through $p$ may be 2-dimensional. It is easy to check that by normality the set
 $\hat P:= \{p\in P: \dim \ell_p = 2\}$ is finite. Restricting to $P\setminus \hat{P}$,
 the dual to $\CF$ is also a family of curves.
\item If $\CF$ is a normal two-dimensional family of curves whose dual is a family
 of curves as well, it is easy to check that the dual is almost normal. However,
 it need not be normal. In the application in Section 2 the family $\CF$ will be
  definable in a strongly minimal structure. We leave it as simple exercise to
  check that in that case replacing $p$ with $p/E$ for an appropriate definable
   equivalence relation $E$ both $\CF$ and its dual are normal. See also Fact \ref{rich}.
\item There is no harm assuming that $\CF$ is only almost normal (whose dual is a family of curves), in which case
 its dual is easily checked to be almost normal as well.

\end{enumerate}

The next definition, which is quite technical, is motivated by the
following argument: Our ultimate goal is to prove that the existence
of an almost normal two-dimensional $\CN$-definable family $\CF$ of
plane curves contradicts the stability of $\CN$. The idea is to
start with a curve $X$, and a curve $C_{q_0}$ from $\CF$
touching $X$ at a generic point, $p_0$. By Lemma \ref{tangency2}, there are
curves arbitrarily close to $C_{q_0}$ which either intersect $X$
nowhere near $p_0$ or intersect it twice there. If there were no
other intersection points of these curve with $X$, this gap in the
number of intersection points with $X$ can be seen to contradict
stability. However, the problem is that the curves near $C_{q_0}$
might intersect $X$ elsewhere and in order to control the total
number of intersection points we need to investigate what happens
near those points. Such points will be called {\em special for
$(q_0,X)$}.

\begin{defnn}\vlabel{sp}{\em  Let $X\sub P$ be a definable curve.
 Given $q\in Q$, a point $p\in M^k$ is called {\em special for
$(q,X)$} if there exist $q'\in \nu_q$, and $p'\in \nu_{p}$ such that
$p'\in X\cap C_{q'}$} (equivalently, for every neighbourhoods $U, V$
of $q,p$, respectively, there exist $q'\in U, p'\in V$, such that
$p'\in X\cap C_{q'}$).

A similar definition can be given even if some of the coordinates of
$p$ are taken to be $\pm\infty$ (and then $\nu_{+\infty},
\nu_{-\infty}$ have an obvious meaning).
\end{defnn}

The point $p$ in the above definition can be seen as an asymptotic
direction of $X$, in the sense of the family $\CF$, and the
curve $C_q$ can be seen as an ``$\CF$-asymptote of $X$ at $p$''.
For example, if $\CF$ is the family of affine lines in $\reals^2$,
given by $C_{a,b}=\{y=ax+b\}$, and $X$ is the curve $y(x-1)=1$ then
the point $p=\la 1,+\infty\ra$ is special for $(\la 0,1\ra,X)$ (that
is, the curve $y=1$ is asymptotic to $X$ at $(1,+\infty)$). Indeed,
curves of the form $y=1+\epsilon$ will intersect $X$ near $+\infty$.
\\

\noindent The following are easy observations concerning special points:
\begin{enumerate}
\item  If $p\in C_q\cap X$ then clearly $p$ is a special point for
$(q,X)$.

\item If $p\in M^k$ is a special point for $(q,X)$ then $p\in \cl(X)\sub \cl(P)$.
However, if $\CF$ is not given continuously, it is possible that $p$
is not in $\cl(C_q)$.

\item
 The set of
special points for $(q,X)$ is clearly definable.

\item If $p$ is a special point for $(q,X)$ then, by curve selection,
 there exist an interval
$(a,b)$ and a  definable curve $\gamma:(a,b)\to Q$ tending to $q$ as
$t$ tends to $a$, such that $p$ is in the limit set $\CF(\gamma)$. In particular,
if $q$ is generic in $Q$ then $p\in
\CF(\gamma)=\cl(C_q)$ (see (2) of the previous section).

\end{enumerate}

The following lemma established the connection between the notions
of a special point and that of a limit set:
\begin{lemma}\vlabel{Special-limit}
Given $\CF$ and a curve $X$, the point $p$ is a special point for
$(q,X)$ if and only if $q$ belongs to the limit set $\CL(X,p)$, of
the dual family $\CL$ along $X$ at $p$ (see (3) of the previous
subsection)
\end{lemma}
\begin{proof} By definition, $p$ is special for $(q,X)$ if and only
if $X\cap \nu_p\cap C_{q'}\neq \0$ for some $q'\in \nu_q$, if and
only if $\ell_{p'}\cap \nu_q \neq\0$ for some $p'\in X\cap \nu_p$.
Using the remark immediately after the definition of limit sets,
this is equivalent to $q$ belonging to the limit set of $\CL$ at
$p$, along the curve $X$.\end{proof}

\begin{lemma}\vlabel{specialpoint}
Let $q_0$ be generic in $Q$ and let
$X\sub P$ be a definable curve not $\CF$-linear near any
of its points. Then:
\begin{enumerate}
\item There are at most finitely many points $p_1,\ldots, p_r$ (including
possibly points with coordinates in $\pm\infty$), which are special
for $(q_0,X)$.
\item There are definable open neighbourhoods $U_1, \ldots, U_r$ of
$p_1,\ldots,p_r$, respectively, and an open neighbourhood $W$ of
$q_0$ such that for every $q\in W$, every intersection point of
$C_q$ with $X$ is contained in one of the $U_i$.
\end{enumerate}
\end{lemma}
\proof (1) As was remarked above, if $p$ is special for $(q_0,X)$
then there is a definable curve $\gamma:(a,b)\to Q$, tending to
$q_0$ as $t$ tends to $a$, such that $p\in
\CF(\gamma)=\cl(C_{q_0})$. Now, if there were infinitely many
special points for $(q_0,X)$ then infinitely many of them were
already in $\Cq\cap X$ and in particular, $\Cq$ would equal $X$ in
some neighbourhood, contradicting our choice of $X$.

(2) If the statement failed then, by curve selection, we would get a
definable curve $\gamma$ tending to $q_0$ with points in
$C_{\gamma(t)}\cap X$ tending away from $p_1,\ldots, p_r$.  By
definition, if $p$ is a limit point for the set $C_{\gamma(t)}\cap
X$ then it is a special point for $(q_0,X)$ (some of the coordinates
of $p$ could be $\{\pm \infty\})$.This would contradict the
assumption that $p_1,\dots,p_r$ were all the special points for
$(q_0,X)$. \qed

The following technical lemma will play a very important role in our
subsequent proof of the main theorem.

\begin{lemma}\vlabel{asymptotes} Let $X\sub P$ be an $A$-definable curve
of rank $k>2$, $p_0$ generic in $X$ and assume that for $q_0\in Q$,
the curve $\Cq$ touches $X$ at $p_0$ (here we identify $P$ locally
with an open set in $M^2$), and in addition,
$\Mdim(p_0,q_0/\emptyset)=3$. If $p$ is a special point for
$(q_0,X)$ then $\Mdim(p/A)=1$. In particular, $p$ is a generic point
of $X$.
\end{lemma}
\proof As was pointed out above, $p\in \cl(X)$ and hence
$\Mdim(p/A)\leq 1$. By \ref{very normal}, we may assume that in a
neighbourhood $U$ of $p_0$ and $W$ of $q_0$, the family $\CF$ is
nice. Consider $\tau(X\cap U)$ and recall that $q_0$ is generic in
$\tau(X\cap U)$ over $A$ (by the Duality Theorem).

Assume, towards a contradiction, that $\Mdim(p/A)=0$. Hence, $q_0$
is still generic in $\tau(X\cap U)$ over $Ap$. As we pointed out
in (3) above, $q_0$ is in the $Ap$-definable limit set $\CL(X,p)$.
But then there is a neighbourhood $W_1\sub W$ of $q_0$ in
$\tau(X\cap U)$, such that every $q$ in $\tau(X\cap U)\cap W_1$
belongs to $\CL(X,p)$. The set $\CL(X,p)$ has  dimension one and
therefore $\tau(X\cap U)$ equals to $\CL(X,p)$ in some open set
$W_2\sub W_1$.

The curve $\CL(X,p)$ is a member of the family
of limit sets $\tilde \CL$, hence by Theorem \ref{limit sets},
its dimension is at most $2$. It follows that the rank of
$\CL(X,p)\cap W_2=\tau(X\cap U)\cap W_2$ is
at most $2$. However, by \ref{rank} (after possibly shrinking
$W_2$), the rank of
$\tau(X\cap U)\cap W_2$ equals that of $X$, contradicting our
assumption that Rank $X\geq 3$.\qed

%
%


\bigskip
\noindent {\bf Remark } The above lemma implies that, under the same
assumptions,  a point in $Fr(X)$ cannot be special for $(q_0,X)$
(see Figure 7). Similarly, it can be shown that such special points
cannot have $\pm\infty$ as one of their coordinates. This can be
seen as saying that a generic tangent curve to $X$ cannot be
asymptotic to $X$ at $\pm \infty$.

\begin{center}
\includegraphics[height=1.2in,width=2in]{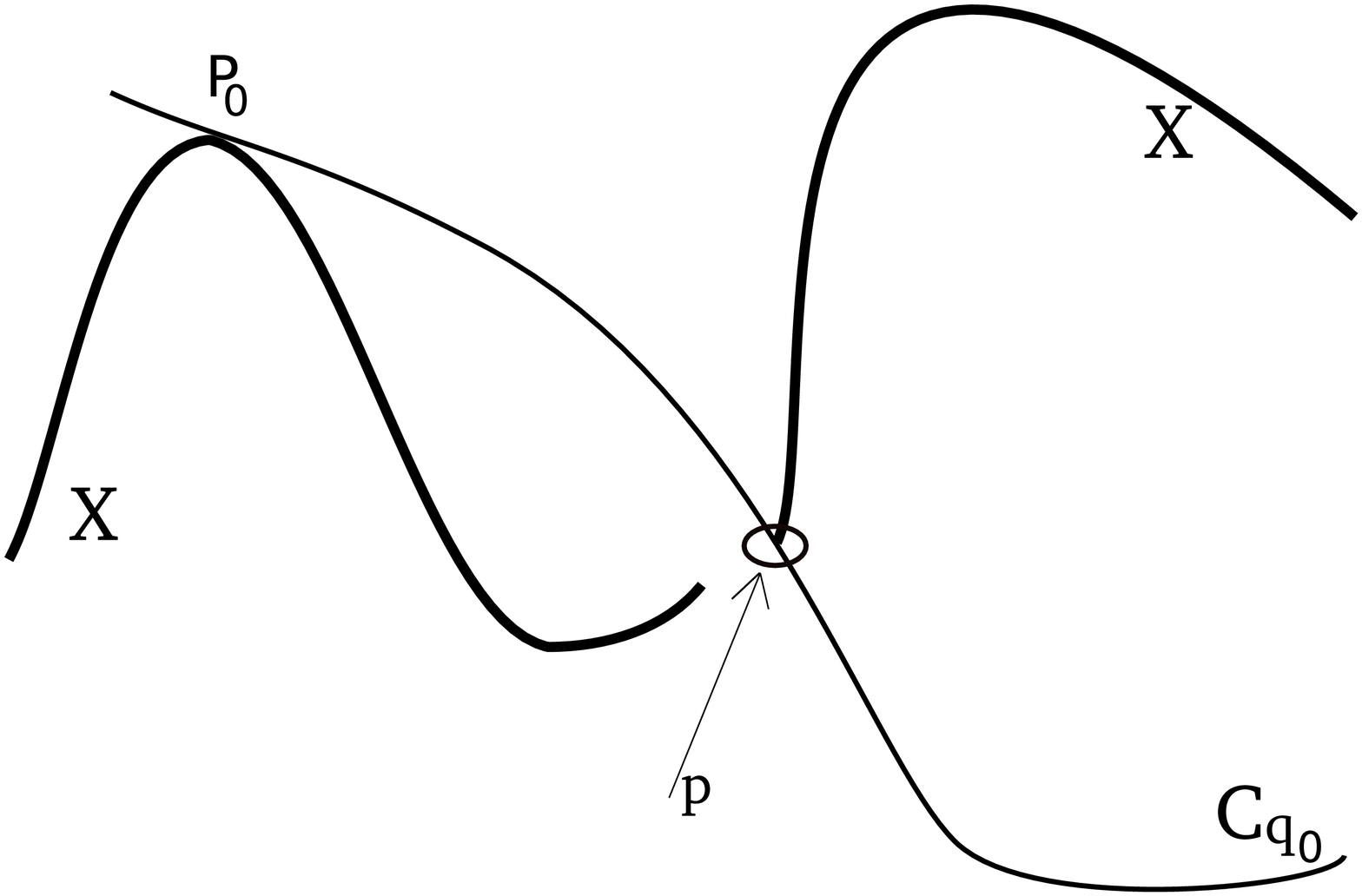}

Figure 7:\\ A problematic special point
\end{center}
We will need the following lemma:

\begin{lemma}\vlabel{transversal2} Let $\CF$, $\CL$, $X$, $p_0$ and
$q_0$ be as in the last lemma. Assume that $p$ is a special point
for $(q_0,X)$ such that $\Mdim(p/p_0,q_0)=1$. Then, for all $q\in
\nu_{q_0}$,  the curve $C_q$ intersects $X$ exactly once in $\nu_p$.
\end{lemma}
\pf  By Lemma \ref{asymptotes}, $p$ is generic in $X$ and by Lemma
\ref{specialcase}, $\C_{q_0}$ does not touch $X$ at $p$.

The fact that $p$ is a special point for $(q_0,X)$ implies (see
(4) after the definition of special points) that $p$ is in
$\cl(\Cq)$. However, since we assume that $\Mdim(p/q_0)=1$, $p$
must belong to $\Cq$. We now apply Lemma \ref{transversal} to
$\CF_{U,W}$ for some neighbourhood $U$ of $p$ and $W$ of $q_0$, such
that $\CF_{U,W}$ is nice.\qed

\section{Stable one dimensional theories in o-minimal structures}\label{SecStable}

We now return to the original setting of the main theorem. Namely,
$\CM$ is an o-minimal structure (expanding a real closed field), $\CN$
is definable in $\CM$ such that $\Mdim(N)=1$. In this section we treat
 the case that $\CN$ is stable. First note that by cell
 decomposition, and the fact that we assume that $\CM$ expands a
 field, we may assume that $N$ is a subset of $M$ (later on, in
 Appendix B, we will show why this can be assumed even without the
 field assumption)

We will need the following facts about stable theories interpretable
 in o-minimal structures:
\begin{fact}
Assume $\CN$ is a stable structure interpretable in $\CM$. Then:\item
\begin{enumerate}
    \item $\CN$ is superstable and its $U$-rank is at most the o-minimal
dimension of $N$ (see \cite{Ga}).
\item If $U(\CN)=1$ and $\CN$ is non trivial then $R^{\infty}(\CN)=1$.
 If in addition $\CN$ is non locally modular then $RM(\CN)=1$
 (see Chapter 2, Corollary 3.3 in \cite{Pi}).
\end{enumerate}
\end{fact}

Our goal in this section is to prove;

\begin{theorem} \vlabel{stable}If $\CN$ is stable then the $acl$-geometry is
necessarily locally modular.
\end{theorem}
\pf By the properties above, $U(\CN)=1$. We assume towards
contradiction that the $acl$-geometry is not locally modular. Hence,
by the same properties, we may assume that $\CN$ is strongly
minimal. Our strategy is to find a strongly minimal field definable in $\CN$.
This leads to an immediate contradiction since on the one hand, by a classical
 result of MacIntyre's any such field must be algebraically closed (see, e.g.
  \cite{Po1}) whereas any strongly minimal field definable in $\CN$ is necessarily
  1-dimensional in the o-minimal sense and hence it is real closed
   (see \cite{vdDr}).

Recall that we use $\Ndim$ to denote the dimension of tuples and
sets in the sense of the geometric structure $\CN$. Since $\CN$ is
strongly minimal this is the same as Morley Rank.

 Non-local modularity of
$\CN$ supplies us with a two-dimensional family of $\CN$-definable
plane curves, to which we intend to apply the results of the
previous section. So we remind:

\begin{defnn}
A definable (possibly in $\CN^{eq}$) {\em pseudoplane} is given by
two sets $P,Q$, together with a definable set $F\sub P\times Q$,
such that the following holds:

 (i) $\mr(P)=\mr(Q)=2$.

 (ii) For every $q\in Q$, $\mr(F(P,q))=1$ (where $F(P,q)=\{p\in P:\la p,q\ra \in F\}$).

(iii) For every $p\in P$, $\mr(F(p,Q))=1$ (where $F(p,Q)=\{q\in Q:\la
p,q\ra \in F\}$).

(iv) For every $q_1\neq q_2$ in $Q$, $F(P,q_1)\cap F(P,q_2)$ is
finite.

(v) For every $p_1\neq p_2$ in $P$, $F(p_1,Q)\cap F(p_2,Q)$ is
finite.
\end{defnn}

Given a pseudoplane as above, $p\in P, q\in Q$, we let $C_q=F(P,q)$
and $\ell_p=F(p,Q)$. We also let $\CF$ denote the family of curves
$\{C_q:q\in Q\}$ in $P$, and let $\CL$ denote the dual family.

The following can be found in \cite{Zil5} (or see Proposition 1.7 on p.155 of \cite{Pi}):
\begin{fact}\label{rich} If $\CN$ is any strongly minimal nonlocally
modular structure, then a pseudoplane is interpretable in $\CN$.
Moreover,  (see, e.g. \cite{Hr5}), $P$ may be taken to be a subset
of $N^2$ (but without EI, we cannot assume at the same time that $Q$
is a subset of $N^2$ as well).
\end{fact}

 From now on we work with a fixed such pseudoplane that we assume is
interpretable in  $\CN$. From the point of view of $\CM$, the sets
$P$ and $Q$ are 2-dimensional definable (using $\CM$-definable choice for $\CN$)
sets, with $\CF$ and its dual $\CL$ definable as well. Clearly,
both $\CF$ and $\CL$ are normal families. We may assume that $P,Q$
and each generic curve in $\CF$ and $\CL$ are of Morley degree
$1$. We may also assume that $P,Q$ and $\CF$ are all $\0$-definable.

The main stage in defining a field in $\CN$ is
to define a group. This will be achieved in two steps.
\\

\noindent {\bf Step I } Constructing an $\CN$-definable curve $X\sub
P$ of rank $k\geq 3$, which has a
sufficiently generic curve from $\CF$ touching it at some generic
$p_0$.
\\

We fix $p_0=\la x_0,y_0\ra  $ an $\CM$-generic point in $P$.
 and $3$ points:
$$p_3=\la z_2,y_0\ra,\,\,  p_2=\la z_2,z_1\ra\,\,p_1=\la x_0,z_1\ra,$$
such that $x_0,y_0,z_1,z_2$ are $\CM$-generic and $\CM$-independent
in $N$, all close to $x_0$ (we will see later how close we need them further to be).

We let $Q_i$, $i=0,1,2,3$, be the set of all $q\in Q$ such that
$p_i\in C_q$. So, each $Q_i$ is $\CN$-definable and has dimension (and Morley rank)
one. We may assume that each $Q_i$ has Morley degree one. For
$q_1\in Q_1, q_2\in Q_2$ and $q_3\in Q_3$ we let
$$C_{q_3 q_2 q_1}=C_{q_3}\circ C_{q_2}^{-1}\circ C_{q_1}.$$ Notice
that if $q_1,q_2,q_3$ is an $\CM$-generic and $\CM$-independent
sequence with $q_i\in Q_i$, then $p_0$ is
$\CM$-generic in $C_{q_3 q_2 q_1}$ (over $q_3 q_2 q_1$). Fix
$X=C_{q_3q_2q_1}$ for such an $\CM$-generic, $\CM$-independent sequence
of $q_i$.

\noindent{\bf Case 1}  The rank of $X$ over
$p_1,p_2,p_3$ is smaller than $3$.

In that case, we obtain an $\CN$-definable infinite field. To do this, recall the
following theorem of Hrushovski's (see, e.g., \cite{Bo}):
\begin{fact}
Suppose that $M$ is a strongly minimal structure and there exists is a collection
 of tuples such that $\mr(g_1,g_2,g_3,b_1,b_2,b_3)=5$ and:
\begin{itemize}

\item $\mr(g_i)=2$ and $\mr(b_i)=1$ for all $1\leq i \leq 3$.
\item The $g_i$ are pairwise independent but $g_3\in \acl(g_1,g_2)$.
\item $\mr(g_1,b_1,b_2)=\mr(g_2,b_2,b_3)=\mr(g_3,b_1,b_3)=3$.
\end{itemize}
Then a strongly minimal field is interpretable in $M$. We call such
 a collection of tuples a {\em field configuration}.

\end{fact}
We will show that under the above assumptions we can find a field
configuration. Note, first, that the rank of $X$ over
  $p_1,p_2,p_3$ is at least 2
(because fixing $q_2,q_3$ we obtain, moving $q_1$, a family of the
same dimension as $Q$). So assume that the rank of $X$ is precisely
2. In particular this is the rank of $C_{q_2}^{-1}\circ C_{q_1}$,
i.e. $\mr (Cb(C_{q_2}^{-1}\circ C_{q_1})) = 2$.  Let $x\in N$ be
generic over $q_1,q_2$ and choose $q_1(x)\in C_{q_1}(x)\setminus
\Nacl(x)$, where $C_{q_1}(x):= \{y:\la x,y\ra\in C_{q_1}\}$. Let
$q_2^{-1}q_1(x)$
 be any point in $(C_{q_2}^{-1}\circ
C_{q_1})(x)\setminus \Nacl(x)$. Then
\[
Q:=(Cb(C_{q_1}), Cb(C_{q_2}), Cb(C_{q_2}\circ C_{q_1}^{-1}), x,
  q_1(x), q_2^{-1}{q_1}(x))
\]
form a field configuration over $p_1,p_2,p_3$. To see
this, note first that $\mr(Q)=5=\mr (q_1,q_2,x)$ since all the other
element in $Q$ are, by construction, algebraic over them The first two bullets in the definition of a field
configuration being immediate from the construction of $Q$ we only
have to check the third condition. But this follows from the fact
that $x$ is generic and (e.g.) $C_{q_1}$ is a generic curve through
$(x,q_1(x))$ whence
$\mr(Cb(C_{q_1})/(x,q_1(x)))=1$. The rest follows in a similar way.
As we have already said, the existence of a strongly minimal (hence
$\CM$-1-dimensional)  field interpretable in $\CN$ leads to a
contradiction with the stability of $CN$.

So we are reduced to:

 \noindent{\bf Case 2} The rank of $X$ over $p_1,p_2,p_3$ is at least $3$.

Let $R$ denote the field underlying the structure $\CM$. Because $p_0, p_1, p_2, p_3 $ are
 generic in $X, C_{q_1},C_{q_2},
C_{q_3}$, respectively, the curves are all $C^1$-functions in
neighbourhoods  of those points. If we take $x_0,y_0, z_1, z_2$ sufficiently close to each
other and $q_1,q_2,q_3$ sufficiently close to each other then by
continuity arguments, the derivative of $X$ at $p_0$ is close to
those of the $C_{q_i}$ at the $p_i$. Therefore, by Lemma \ref{diff2} and continuity,
we can find $q_0\in Q_0$ such that $\Cq$ has the same derivative as
$X$ at the point $p_0$. By Theorem \ref{diff3}, the curve $\Cq$
touches $X$ at $p_0$, and $\Mdim(p_0,q_0/\emptyset)=3$. {\bf End of
Step I}
\\

\noindent{\bf Step II} Counting Intersection points of $X$ with $\CF$-curves.

\bigskip
We now fix the curve $X$, $p_0$ generic in $X$ and $C_{q_0}$
touching $X$ at $p_0$ as given by Step I. We will investigate  the
intersection with $X$ of $\CF$-curves close (in $Q$) to $q_0$.
 By Lemma \ref{specialpoint}, all intersection points of such curves with $X$
 appear near the special points for $(q_0,X)$. By the first part of the same
  lemma, there are at most finitely many such special points
  $s_0=p_0, s_1,\ldots, s_r$

\bigskip
\noindent{\bf The General Case:} All $s_i$, except $s_0$, have
$\Mdim(s_i/p_0,q_0)=1$.

By Lemma \ref{specialpoint} there exist neighbourhoods $W$ of $q_0$
and pairwise-disjoint $U_i$ of $s_i$, $i=0,\ldots, r$, respectively,
such that for every $q\in W$, every intersection point of $C_q$ with
$X$ is in one of the $U_i$. Furthermore, by Lemma \ref{transversal2}
we can choose $W$ and the $U_i$ in such a way that for each $q\in W$
and $i\neq 0$, the curve $C_q$ intersects $X$ in $U_i$ exactly once.

Now, by Lemma \ref{tangency2}, we can find $q_1, q_2\in W$ both
$\CN$-generic over all the data, such that, after possibly shrinking
$U_0$ further, we have  $C_{q_1}\cap X\cap U_0=\emptyset$ while
$|C_{q_2}\cap X\cap U_0|>1$.

 It follows  that
$$|X\cap C_{q_1}|=r,$$ while
$$|X\cap C_{q_2}|\geq r+2.$$

This contradicts the assumption that the Morley degree of $Q$ is $1$.

\bigskip


\bigskip

\noindent{\bf The Special Case:} There exists $s_i$, $i=1,ldots,r$
such that $\Mdim(s_i/p_0q_0)=0$.  The crucial part of the proof is:

\begin{lemma}\vlabel{degenerate}  If, under the above assumptions, $\Mdim(s_1/p_0,q_0)=0$
then a group is interpretable in $\CN$.
\end{lemma}
\proof We still work with the $4$ points
$$p_0=\la x_0,y_0\ra \, , \, p_1=\la x_0,z_1\ra\, , \, p_2
=\la z_2,z_1\ra\, , \,p_3=\la z_2,y_0\ra,$$
such that $x_0,y_0,z_1,z_2$ are $\CM$-generic and $\CM$-independent.
As above, we let $Q_i$, $i=1,2,3$, be the set of all $q$.
such that $p_i\in C_q$. So, each $Q_i$ is $\CN$-definable and has
dimension one. We denote elements of $Q_2$ by $r, r_1,r_2$ etc.

For $q_1\in Q_1, r\in Q_2$ and $q_3\in Q_3$ we use
$C_{q_3rq_1}=C_{q_3}\circ C^{-1}_r\circ C_{q_1}$ as before.  We
denote the composition curves $C_{q_3}\circ C_{r_2}^{-1}$ and
$C_{r_2}^{-1}\circ C_{q_1}$ by $C_{q_3 r_2}$ and $C_{r_2 q_1}$,
respectively.

We fix $\CM$-independent $\hat{q_3}, \hat{r}, \hat q_1 $
$\CM$-generic in $Q_1,Q_2, Q_3$, respectively, and, as before,
$X=C_{\hat q_3 \hat r \hat q_1}$, and $\Cq$ touches $X$ at $p_0$,
 with $\Mdim(p_0,q_0/\emptyset)=3$. We assume also that there is
$s_1\in P$ special for $(q_0,X)$ such that $s_1\in \Mdcl(p_0,q_0)$.
By Lemma \ref{asymptotes}, $s_1$ is generic in $X$ over the
parameters defining $X$ (so $s_1\notin Fr(X)$).

Our goal now is to construct a group configuration in $\CN$. First,
since $$\Mdim(\hat q_1,\hat r,\hat q_3/p_0,p_1,p_2,p_3)=3,$$ there
are relatively open neighbourhoods $Q_1'\sub Q_1\,\, ,Q_2'\sub Q_2$
and $Q_3'\sub Q_3$  of $\hat q_1,\hat r,\hat q_3$ respectively
(definable in $\CM$!!) such that:
\begin{align*}
\tag{*}
&\mbox{for every\,} q_1\in Q_1', r\in Q_2',q_3\in Q_3',\\
&\mbox{ if $\Cq$ touches $C_{q_3 r q_1}$ at $p_0$ then $s_1\in
C_{q_3 r q_1}$}
\end{align*}
 This can be done because touching
is an $\CM$-definable notion and there is a $p_0$-definable function
sending $q_0$ to $s_1$.
\\

We first claim that for $r_1, r_2\in Q_2$ which are
$\CM$-independent and $\CM$-generic (over all mentioned parameters),
the curve $C_{r_1r_2}:=C_{r_1}C_{r_2}^{-1}$ (we will use this notation whenever
$r_1,r_2\in Q_2$) is $C^1$ in a neighbourhood
of $\la z_1,z_1\ra$: Indeed, because
$\la z_2,z_1\ra$ is generic in both $C_{r_1}, C_{r_2}$, it
is $C^1$ in both. It is therefore sufficient to see that
that the germ of $C_{r_1r_2}$ at $\la z_1,z_1\ra$ equals
the composition of the germ of $C_{r_1}$ at $\la z_2, z_1\ra$ with
the germ of $C_{r_2}^{-1}$ at $\la z_1,z_2\ra$. For that, we need to
verify that there is no $z_2'\neq z_2$ such that $\la z_2',z_1\ra
\in \cl(C_{r_1})\cap \cl(C_{r_2})$.

If such $\la z_2',z_1\ra $ were in $Fr(C_{r_i})$ then $z_1$ would
not be generic over the parameters defining $C_{r_1}$,
contradiction. So we may assume that $\la z_2',z_1\ra \in
C_{r_1}\cap C_{r_2}$, with $z_2'\neq z_2$. Since the set $\{z:\la
z,z_1 \ra\in C_{r_1}\}$ is finite we have $z_2,z_2'\in
\Mdcl(z_1,r_1)$. Because $r_2$ was chosen to be generic in $Q_2$
over $z_1,r_1,p_2$, we get infinitely many $r\in Q_2$ such that $\la
z_2,z_1\ra, \la z_2',z_1\ra\in C_r$. This contradicts the normality
of $\CL$.
\\

We can now describe the idea for the rest of the proof.
\\

 \noindent {\bf The idea: } Generically, we will be able
to recognise in $\CN$ when $C_{r_1r_2}$ and
$C_{r_3r_4}$ have the same derivative at $\la z_1,z_1\ra$,
for $r_1,r_2,r_3,r_4\in Q_2$. This is done via the following
relation:

 {\em $C_{r_1r_2}$ and $C_{r_3r_4}$ are
tangent to each other at $\la z_1,z_1\ra$ iff there are $q_1,q_1'\in
Q_1'$ and $q_3,q_3'\in Q_3'$ such that the curves $C_{q_3 r_1 q_1}$,
$C_{q_3 r_2 q_1'}$, $C_{q_3' r_3 q_1}$ and $C_{q_3' r_4 q_1'}$ all
touch $\Cq$ at $p_0$.}

Indeed, notice that if $C_{q_3 r_1 q_1}$ and $C_{q_3' r_3 q_1}$
touch $\Cq$ at $p_0=\la x_0,y_0\ra$ then in particular they have the
same derivative at $p_0$ and therefore $C_{q_3 r_1}$ and $C_{q_3'
r_3}$ also have the same derivative at $\la z_1,y_0\ra$. Similarly,
$C_{q_3 r_2}$ and $C_{q_3' r_4}$ have the same derivative at $\la
z_1,y_0\ra$. By composition, it now follows that $C_{r_1r_2}$ and
$C_{r_3r_4}$ have the same derivative at $\la z_1,z_1\ra$. Why is
the above relation $\CN$-definable ? By $(*)$, for independent
generic $q_1,r,q_3$ in $Q_1',Q_2',Q_3'$, respectively, tangency of
$C_{q_3rq_1}$ to $\Cq$ at $p_0$ can be detected by the fact that
$s_1$ belongs to $C_{q_3 r q_1}$.

Finally, using the fact that there is an underlying real closed
field the above tangency relation will allow us to obtain a group
configuration in $\CN$.
\\

\noindent{\bf The actual argument} : For $r_1,r_2,r_3,r_4\in Q_2$ we
define the relation

\begin{align*}
&T(r_1,r_2,r_3.r_4)\Leftrightarrow \\
&\mbox{\em there are  $q_1,q_1'\in Q_1$  and $q_3,q_3'\in Q_3$}
\mbox{\,\,such that}\\
& s_1\in C_{q_3 r_1 q_1}\cap C_{q_3 r_2 q_1'}\cap C_{q_3' r_3q_1}
\cap C_{q_3' r_4 q_1'}
\end{align*}

$T$ is clearly definable in $\CN$. The main properties of $T$ are:

\begin{claim}\vlabel{propT}
Denote $A:=\{s_1,p_0, p_1,p_2,p_3,q_0\}$ then:
\begin{enumerate}
\item For $\la r_1,r_2\ra$ and  $\la r_3,r_4\ra$, each $\CM$-generic  in $Q_2'\times Q_2'$
and sufficiently close to $\hat r$, if the curves
$C_{r_1r_2}$ and $C_{r_3r_4}$ have the same
derivative at $\la z_1,z_1\ra$ then $\la r_1,r_2,r_3,r_4\ra\in T$.
\item If $\la r_1,r_2,r_3\ra$ is an $\CN$-generic triple of elements from $Q_2$
then there exists $r_4\in Q_2$, $r_4\in \Nacl(r_1r_2r_3A)$,  such
that $\la r_1,r_2,r_3,r_4\ra \in T$.
\end{enumerate}
\end{claim}

The proof of (1): Assume that $C_{r_1r_2}$ and
$C_{r_3r_4}$ have the same derivative at $\la z_1,z_1\ra$
(and the $r_i$ are close to each other).

Pick $q_3\in Q_3$ to be $\CM$-generic over $r_1,r_2,r_3,r_4, q_0$
and close to $\hat q_3$. Because $r_1$, $r_3$ are close to each
other, the derivatives of $C_{q_3 r_1}$ and $C_{q_3 r_2}$ at $\la
z_1,y_0\ra$ are close to that of $C_{\hat q_3\hat r}$. We can now
find $q_1, q_1'\in Q_1'$ such that $C_{q_3 r_1 q_1}$ and $C_{q_3 r_2
q_1'}$ have the same derivative at $p_0$ as $X=C_{\hat q_3\hat r\hat
q_1}$ does and hence also $\Cq$ (this is possible by Lemma \ref{diff2} and
 continuity arguments in $\CM$). By Lemma \ref{diff3}, the curve
$\Cq$ touches $C_{q_3 r_1 q_1}$ and $C_{q_3 r_2 q_1'}$ at $p_0$. It
follows from $(*)$ that $s_1$ belongs to both $C_{q_3 r_1 q_1}$ and
$C_{q_3 r_2 q_1'}$. It is not hard to see that $q_1$ and $q_1'$ are
each $\CM$-generic over $r_1,r_2,r_3,r_4,q_0$ (because each of them
is inter-definable with $q_3$ over the $r_i$ and $q_0$). For the
same reason,
 we can find $q_3'\in Q_1'$ such that the derivative of
$C_{q_3' r_3 q_1}$ at $p_0$ equals to that of $\Cq$.

Also, by our assumption on the $r_i$, the curves
$$C_{q_1}^{-1}C_{r_1}C_{r_2}^{-1}C_{q_1'}\,\,\,\mbox{ and }
\,\,\,
 C_{q_1}^{-1}C_{r_3}C_{r_4}^{-1}C_{q_1'}$$ have the same
 derivative at $p_0$ and therefore

$$C_{q_1}^{-1}C_{r_1}C_{q_3}^{-1}C_{q_3}C_{r_2}^{-1}C_{q_1'}\,\, \mbox{ and }
C_{q_1}^{-1}C_{r_3}C_{q_3'}^{-1}C_{q_3'}C_{r_4}^{-1}C_{q_1'}$$ have
the same derivative at $\la x_0,x_0\ra$.

Because $C_{q_3 r_1 q_1}$, $C_{q_3 r_2 q_1'}$ and $C_{q_3' r_3 q_1}$
all have the same derivative at $p_0$ it follows from the above that
$C_{q_3' r_4 q_1'}$ has the same derivative as well and therefore,
by Lemma \ref{diff3}, it touches $\Cq$ at $p_0$. Hence $s_1\in
C_{q_3' r_4 q_1'}$ and so $\la r_1,r_2,r_3,r_4\ra\in T$.

The proof of (2). We prove: for every generic tuple $r_1,r_2,r_3\in
Q_2^3$ there are finitely many (and at least one) $r_4$ such that
$\la r_1,r_2,r_3,r_4\ra \in T$. By the strong minimality of $\CN$
this is a first order statement. Because $Q_2\times Q_2\times Q_2$
has Morley degree $1$, it is sufficient to prove this for any
particular $\CN$-generic triple in $Q_2^3$. This can be done
similarly to (1). Fix $\la r_1,r_2,r_3\ra\in Q_2'\times Q_2'\times
Q_2'$, $\CM$-generic over $A$, where the $r_i$ are close to $\hat
r$. As in the proof of (1), there is $r_4\in Q_2'$ such that
$C_{r_1r_2}$ and $C_{r_3r_4}$ have the same derivative at $\la z_1,
z_1\ra$. Moreover, each pair $r_1,r_2$ and $r_3, r_4$ is generic, so
by (1) we have $\la r_1,r_2,r_3,r_4\ra \in T$. Actually, as the
argument in (1) shows, to witness this we may choose $q_3\in Q_3'$
as
generic as we wish, and then find the appropriate $q_1,q_1',q_3'$.  \\

To finish the proof of (2), we need to show that $r_4\in
\Nacl(r_1,r_2,r_3,A)$. We first prove:

 \noindent {\bf Claim } $r_4\in \Nacl(r_1 r_2
r_3 q_3A)$.
\\

Indeed, once we choose $q_3$, the family $Y=\{C_{q_3 r_1 q}: q\in
Q_1\}$ is one-dimensional, through $p_0$. Notice that $s_1\notin
\Mdcl(p_0)$ for otherwise there will be infinitely many curves from
$\CF$ going through $p_0$ and $s_1$, contradicting the normality of
$\CL$.
 Since $q_3, r_1$ were
chosen generically over $A$  there are at most finitely many curves
in the family $Y$ which go through both $p_0$ and $s_1$. Hence,
$q_1\in \Nacl(r_1q_3A)$. Similarly, $q_1'\in \Nacl(r_2 q_3 A)$,
$q_3'\in \Nacl(r_3q_1)$ and hence $q_3'$ is also in $\Nacl(r_1r_2r_3
q_3 A)$. Finally, we may conclude in a a similar fashion that
$r_4\in \Nacl(r_1r_2r_3 q_3 A)$, as needed.
\\

Assume towards a contradiction that $r_4\notin \Nacl(r_1 r_2 r_3
A)$. Then, by the Exchange Principle and the claim above, $q_3\in
\Nacl(r_1 r_2 r_3 r_4A)$ and in particular, $q_3\in \Macl(r_1 r_2
r_3 r_4 A)$. This contradicts our freedom, as discussed above, to
choose $q_3$ from $Q_3'$ to be $\CM$-generic over $r_1r_2 r_3 r_4
A$. End of (2).
\\

This finishes the proof of Claim \ref{propT}, and we can now return to the task
of finding a group configuration in $\CN$. \\

For $r\in Q_2'$ let $d(r)$ denote the derivative of $C_r$ at $\la
z_1,z_1\ra$. Notice that  the derivative of $C_{r_1r_2}$ at
$z_1$ is $d(r_1)/d(r_2)$.
  It follows from the properties of $T$ described in \ref{propT} that, for $r_1,r_2,r_3$
  sufficiently close to each other in $Q_2$ and $\CM$-independent,
  there exists $r_4$ such that $d(r_1)/d(r_2)=d(r_3)/d(r_4)$ and
  $r_4\in \Nacl(r_1r_2r_3A)$.

It is now easy to get a group configuration in $\CN$ using the
(local) group $\{d(r):r\in Q_2'\}$. End of Lemma
\ref{degenerate}.\qed

\bigskip

\noindent {\bf Remark} The assumptions of Lemma \ref{degenerate} are
not vacuous. Assume for example that at every $q\in Q$, every
possible derivative is represented by two different $\ell_p$
through $q$. Given any definable curve $Y\sub Q$, consider the curve
$X=\tau(Y)\sub M^2$. Because every $q\in Y$ has two different
$\CL$-tangent curves $\ell_p$ and $\ell_{p'}$, the duality Theorem
implies that $C_q$ is tangent to $X$ at both $p$ and $p'$. We thus
have $p'$ special for $(q,X)$ and $p'\in \Macl(p,q)$.
\\

We now return to The Special Case, with $X,p_0$ and $q_0$ as before.
Namely, one of the special points for $(q_0,X)$ is in the
$\CM$-definable closure of $p_0,q_0$ (over $p_1,p_2,p_3$). Then, by
Lemma \ref{degenerate} we can interpret a 1-dimensional group in
$\CN$. Because of strong minimality, the group can be assumed to
live on $\CN$, at least locally. Namely, we may assume now that we
have a local 1-dimensional abelian group operation, denoted $+$,
living on $N$.

We now start the process all over again and consider the two
families of curves obtained by compositions
$C_{q_3}C_{q_2}^{-1}C_{q_1}$ and addition (using the group
structure)  $C_{q_3}-C_{q_2}+C_{q_3}$ of the original $C_q$. We
may assume that we are in the case of Lemma \ref{degenerate} for
both families of curves, and then we proceed by defining tangency at
$\la z_1, z_1\ra$  for curves of the form $C_{r_1}C_{r_2}^{-1}$ and
of the form $C_{r_1}-C_{r_2}$ (subtraction in the sense of the group
we just defined). Using the underlying real closed field structure
we can now obtain a definable 1-dimensional field in $\CN$,
contradicting strong minimality.\qed

\section{ The Trichotomy Theorem}\label{Tri}


Let $\CN$ be a structure definable in an o-minimal structure $\CM$
(with a fixed interpretation of $\CN$ in $\CM$). A type $p$ in the
structure $\CN$ is called {\em one dimensional} if it contains an
$\CN$-formula whose $\CM$-dimension is one.

Notice that the notion of dimension depends on the particular
interpretation of $\CN$ in $\CM$. E.g., if $\CN$ is a dense linear
ordering then it can be interpreted in $\CM$ as a subset of $M$, or
as a subset of $M^2$, with the lexicographic ordering. In the first
case, every nonalgebraic 1-type in $\CN$ is one dimensional while in
the second case, there are two dimensional types, but if one fixes
two points in the same fibre and then the interval between them
gives rise to one-dimensional types.

If $p$ is a one-dimensional type in $\CN$ then the restriction of
$acl_{\CN}(\cdot)$ to the realization of $p$ forms a pre-geometry
(Exchange holds because it is true in $\CM$).

We will need the following from \cite{OnPe}:

\begin{defnn}
 Let $\CN$ be a dependent theory (i.e. a theory with NIP). A definable set
  $X$ is \emph{unstable in $\CN$} if there exists an $\CN$-formula $\varphi(x,y)$
  defined over $N$ and $\la a_i, b_i \ra _{i\in \omega}$, $a_i,b_i\in X$,
  such that $N\models
   \varphi(a_i,b_j)$ iff $i < j$. If there is no such formula $\varphi$ then
   $X$ is said to be a \emph{stable in $\CN$}.
\end{defnn}

By Corollary 2.6 of \cite{OnPe},  if $X$ is stable set in a
structure with NIP then it is stably embedded, namely every
$N$-definable subset of $X^k$ is definable with parameters from $X$.

Our main tools in the proof will be the Trichotomy Theorem for o-minimal structures
\cite{PeSt2} and the the following result from \cite{HaOn}:

\begin{theorem}\label{ho}
Let $\CN$ be a 1-dimensional structure definable in an o-minimal structure $\CM$.
For any unstable $X\subseteq N$ there exists $\CN$-definable $X_0\subseteq N$ and
 equivalence relation $E$ with finite classes such that $X_0/E$ with all its
 induced $\CN$-structure is o-minimal. In particular $X_0/E$ is linearly ordered.
\end{theorem}

With this in hand we can now prove:

\begin{introtheorem}\vlabel{1types}
Assume that $\CN$ is a definable structure in an o-minimal $\CM$ and
that $p$ is a complete 1-$\CM$-dimensional $\CN$-type over a model
$\CN_0\sub \CN$. Then:
\begin{enumerate}
\item $p$ is trivial (with respect to $acl_{\CN}$).

Or, there exists an $\CN$-definable equivalence relation $E$ with
finite classes such that one of the following holds:

\item $p$ is linear, in which case either
\begin{enumerate}
 \item $p/E$ is a generic type of a strongly minimal $\CN$-definable
one-dimensional group $G$, and the structure which $\CN$ induces on
$G$ is locally modular.  In particular, $p$ is strongly minimal. Or,
 \item for every $a\models p$ there exists
 an $\CN$-definable ordered group-interval $I$ containing $a/E$.
 The structure which $\CN$ induces on $I$ is a reduct of an ordered
vector space over an ordered division ring.
\end{enumerate}
\item $p$ is rich: For every $a\models p$ there exists
 an $\CN$-definable real closed field $R$ containing $a$
 and the structure which $\CN$ induces on $R$ is o-minimal.
\end{enumerate}
\end{introtheorem}

\begin{proof} Because $p$ is one-dimensional it is contained in an
$\CN$-definable one dimensional set $X$.

\noindent{\bf Case 1}:   $X$ is a stable set in $\CN$.

 It follows
(see earlier discussion) that $X$ is stably embedded in $\CN$. In
this case, we may replace $\CN$ by $X$ and assume that $\CN$ is
one-dimensional and $U(\CN)=1$. By Theorem \ref{stable}, $\CN$ is
necessarily 1-based. It now
follows from \cite{Hr4}
  that either $p$ is
trivial or $p/E$ is the generic type of an $\CN$-definable
locally-modular group $G$, for $E$ an $\CN$-definable finite
equivalence relation.

 Because $p$ is one-dimensional, we must have $\dim_{\CM}(G)=1$
  ($p$ is generic in $G$,
 so every definable subset of $G$ containing $p$ is group-generic
 in $G$, namely, finitely many
 group-translates of
 it cover $G$. Since $\dim p =1$ we can find such a set which
  is 1-dimensional). It is left to show that $G$ is strongly minimal.
 Indeed, if not then some infinite co-infinite set $Y\sub G$ is
 $\CN$-definable. By translating $Y$ to $0\in G$, we may assume that
 $0$ lies on the boundaries of both an infinite component of $Y$ and an infinite component of $G\setminus Y$.
 It is not difficult to see now that the relation $x-y\in Y$ defines an
 ordering on a small neighbourhood of $0$, contradicting stability.

\noindent {\bf Case 2}: $p$ is nontrivial and every $\CN$-definable
one-dimensional set containing $p$ is unstable.

By the non-triviality of $p$,  we can then find an $\CN$-definable
ternary relation $R\sub X^3$ (definable possibly over parameters
outside of $N_0$), witnessing the non triviality of $p$, such that
for every independent $a,b\models p$ there are finitely many $c$
such that $\models R(a,b,c)$, and at least one of those $c$ is in
$p$. Moreover, we may assume that the projection map of $R$ on any
of its two coordinates is finite-to-one, and that the image of this
projection has dimension $2$. By compactness, we may shrink $X$
further so that the projection of $R$ onto the first two coordinates
equals $X^2$ modulo possibly a definable set of dimension one.

Because $X$ is unstable we may apply  Theorem \ref{ho} and find an
$\CN$-definable infinite set $X_0\sub X$ and an $\CN$-definable
finite equivalence relation $E$ on $X_0$ such that the structure
which $\CN$ induces on $X_0/E$ is o-minimal. In order to consider it
as a structure on its own right we need stable embeddedness.

By Lemma 2.3 in \cite{PeSt2}, every closed interval in an o-minimal
structure is stably embedded. However, the very same proof there
gives a stronger result: If $\CM'$ is an o-minimal structure
interpretable in a geometric structure $\CN$, with $\Ndim(M')=1$,
then every closed interval in $\CM'$ is stably embedded in $\CN$.

 By replacing $X_0$ by a smaller $\CN$-definable set we may
 assume that $X_0/E$ is a closed interval in an $\CN$-definable o-minimal structure,
 hence stably embedded in $\CN$.
Let $\CX_0$ denote the induced structure on $X_0/E$.

Our goal now is to transfer the o-minimal structure from $X_0/E$ to
elements realizing $p$ (note that a-priori we do not know that $X_0$
contains $p$). We start with $a\models p$ and $b\in X_0$ independent
from $a$. By our assumptions on $R$, there exists $c$ such that
$R(a,b,c)$ and any two of $\{a,b,c\}$ are inter-algebraic over the
third one. We fix $c$ and consider the relation $\bar R=R(x,y,c)$.

The
  projection of $\bar R$
on the first and second coordinate is finite-to-one and $\bar R$
defines a finite-to-finite correspondence between $p'=tp(a/c)$ and
$tp(b/c)$. The map $\alpha$ which sends $\la x,y\ra\in \bar R$ to
$y/E$ is finite-to-one from $\bar R$ onto an infinite definable
subset of $\CX_0$. Let $\beta$ be the projection map of $\bar R$
onto the $x$-coordinate and for $x\in \beta\alpha^{-1}(X_0)$, let
$\gamma(x)=\min_{\CX_0}(\alpha \beta^{-1}(x))$ (the set on the right
is finite hence it has a minimum in the structure $\CX_0$).

$\gamma$ is  a finite-to-one map from an infinite definable set $Y$
containing $p'$ onto an infinite definable set $Z$ (which we may
assume to be a closed interval) in $\CX_0$. We thus obtain a
definable bijection between $Y/\sim$, for an $\CN$-definable finite
equivalence relation $\sim$, and an infinite closed interval in
$\CX_0$. This bijection allows us to transfer the structure which
$\CN$ induces on $Z$ onto $Y/\sim$, hence the structure which $\CN$
induces on $Y/\sim$ is o-minimal (and stably embedded as well), call
it $\mathcal Y$. The type $p'/\sim$ is generic in $\mathcal Y$.

By the Trichotomy Theorem for o-minimal structures, we may assume
(after possibly shrinking $Y$ further) that the structure of $\CY$
is either trivial, an ordered reduct of a group-interval in an
ordered vector-space, or an o-minimal expansion of a real closed
field.

The non-triviality of $p$ easily implies the non-triviality of $p'$
and therefore of $\CY$: Fix $a_1,a_2\models p'$ independent and $d$
such that $\la a_1,a_2,d\ra \in R$, and define
$$\la x_1,x_2,x_3\ra\in R'\Leftrightarrow \exists y \la x_1,x_2,y\ra
\in R \& \la x_3,d,y\ra\in R.$$ The relation $R'$ witnesses the
non-triviality of $p'$.

If $p$ is linear then, by definition, so is $p'$ and therefore $\CY$
cannot define a real closed field (the pull-back of the family of
affine lines will contradict linearity) hence $\CY$ is an ordered
group interval. If $p$ is rich then there is a definable almost
normal family of curves such that $\la a,a\ra$ belongs to infinitely
many of them (see the discussion following the definition of a rich
type). In that case $\CY$ must be an o-minimal expansion of a real
closed field.
\end{proof}

\bigskip

\noindent We conclude by pointing out that the statement of the theorem in the unstable
case cannot be strengthened so that the whole of $p$ is contained in
a single  o-minimal structure:
 Consider $N$ the expansion of the (unordered) group of the reals
by a predicate for the unit interval. In this theory the ordering is
definable on every interval of finite length. Take the type
 $q=\{x>r: r\in \reals\}$. It is not contained in any definable
ordering, but after fixing $b\models q$ the set $b<x<b+1$ is contained in a definable
 o-minimal structure.

\section{Appendix A: A proof of  v. d. Dries' Theorem on limit sets}

In this section we give we give a proof of Theorem \ref{limit sets}.
 The proof is due to v.d. Dries, and relies in parts on \cite{vdDr1},
 with minor changes intended to make it more self-contained.

We restate the theorem:

\begin{theorem}
 Let $\CM$ be an o-minimal expansion of a field, $\CS$ a definable family
  of definable closed subsets of $M^k$, parameterised by $Z$. Let
  $\tilde \CS:=\{\CS(\gamma)\mid \gamma:(0,1)\to Z, \,\mbox{definable}\}$.
   Then $\tilde \CS$ is $\CM$-definable and $\dim (\tilde \CS \setminus \CS )< \dim \CS$.
\end{theorem}

Because of the presence of field, we may assume that all curves in
$Z$ are given by $\gamma:(0,+\infty)\to Z$. It is not hard to check,
using definable choice, that the collection
 $\{\CS(\gamma)\mid \gamma:(0,+\infty)\to Z\}$ depends only on the class
  $\{\CS_a:a\in X\}$ and not on the parameterisation. Hence, we
  may assume that every set in $\CS$ appears exactly once in
   the family and in particular $\dim \CS = \dim Z$ (otherwise re-parameterise, see \S 3
   of \cite{vdDr1}).

The proof goes through the theory of tame elementary pairs (see \S 8
of \cite{vdDr1}
 for the details). Recall that for an o-minimal expansion $\CN$ of a real closed field, a
  tame extension $\CN \precneqq \CR$ is one in which $\CN$ is Dedekind complete.
   The theory of tame elementary pairs is the theory of those structures
   $(\CR, \CN, \st)$ where $\CR$ is a tame extension of $\CN$ and $\st: R\to N$
    is the standard part map. Note that $\st$ is defined on $V$, the convex hull
    of $N$ in $R$. To simplify the notation, whenever $Y\subseteq R$ we will
     write $\st (Y)$ for $\st(Y\cap V)$.

The crucial property of the theory of tame pairs is:

\begin{proposition}\vlabel{QE}
 Suppose that $(\CR, \CN, \st)\models T_{tame}$ and $Y\subseteq N^k$ is
  definable in $(\CR, \CN, \st)$, then $Y$ is definable in $\CN$. Moreover,
   the theory $T_{tame}$ is complete modulo $Th(\CN)$.
\end{proposition}

From this we readily get a weaker version of the theorem, one which
actually suffices for our purposes:

\begin{lemma}\vlabel{easycase}
 Let $\CM$ be an o-minimal expansion of a field, $\CS$ a $\0$-definable family of
  definable subsets of $M^k$, parameterised by $Z$.
  Let $\tilde \CS:=\{\CS(\gamma): \gamma:(0,+\infty)\to Z, \,\mbox{definable}\}$.
  Then $\tilde \CS$ is $\CM$-definable and $\dim \tilde \CS = \dim \CS$.
\end{lemma}
\begin{proof} By replacing $\CS$ with the family $\{cl(A):A\in
\CS\}$ we may assume that all the sets in $\CS$ are closed (it is
easy to check that the family of limit sets does not change).
 We may assume that $\CM$ is saturated enough to assure that every $M$-definable
  set has (in $M$) a generic point. Let $\tau > M$ and $N = M\la \tau \ra$, the
   prime model over $M\cup \{\tau\}$. Then $M$ is Dedekind complete in $N$ and
   $(\CN,\CM,\st)\models T_{tame}$. Let $\gamma:(0,+\infty)\to Z$ be any curve definable
    over $M$. It is not hard to verify that $\CS(\gamma) = \st \CS_{\gamma(\tau)}$.
     By Proposition \ref{QE}, and by compactness,
      the definable family $\{\st \CS_a: a\in Z(N)\}$ is
       $\CM$-definable. Since every $a\in Z(N)$ is of the form $\gamma(\tau)$ for some
      $\CM$-definable curve $\gamma$, the definability of $\tilde \CS$ in $\CM$ follows.
Because every automorphism of $\CM$ leaves $\CS$ invariant and
therefore also $\tilde \CS$, it follows that $\tilde \CS$ is
actually $\0$-definable.

Now let $Z_1$ be a parameter set for $\tilde \CS$ (so $Z_1$ is
$\0$-definable in $\CM$). Again, we assume that every set in $\tilde
\CS$ is represented exactly once in the family.
 Let $b\in Z_1(M)$ be generic over $\0$. There is some $\gamma:(0,+\infty)\to Z$ such that
 $\tilde \CS_{b} = \CS(\gamma) = \st \CS_{\gamma(\tau)}$.
If $\gamma(\tau)\in M^k$ then $\gamma$ is eventually constant
and hence $\CS_{\gamma}\in \CS$. Assume then that
$\gamma(\tau)\notin M$ and let $\mathcal P = \la \gamma(\tau) \ra$,
the model generated by $\gamma(\tau)$ over $\0$. Consider the tame
pair $(\mathcal P, \st \mathcal P, \st)$.
 It follows from the
completeness part of
  Proposition \ref{QE}
  that $\st \CS_{\gamma(\tau)}=\tilde \CS_{b'}$ for
  some $b'\in Z_1(\st \mathcal P)$. By our assumption, we must have $b'=b$.

Since  $\gamma(\tau) \in Z$, we have  $\dim(\gamma(\tau)/\0) \leq
\dim(\CS)$, and therefore
$$\dim(\CS)\geq \dim(\mathcal P/\0)\geq
\dim(\st \mathcal P/\0).$$

We conclude that $\dim(b/\0)\leq \dim(\CS)$  therefore
$\dim(\tilde(S))=\dim(Z_1)\leq \dim(\CS)$.\end{proof}

As mentioned above, the previous lemma suffices for our purposes.
For the sake of
 completeness, however, we give the proof of v.d. Dries' theorem in its full strength.

The main part of the proof is the following:

\begin{lemma} For $\CS$ as above, assume that all sets in $\CS$ are
closed and $\0\notin \CS$. Let $M\prec N\prec R$, be such that $M,N$
are Dedekind complete in $R$, and let $\st: R\to N$ be the
associated standard part map, $V$ the convex hull of $N$ in $R$. Let
$a\in Z(R)\cap V^n$ be such that $M\la a \ra\subseteq V$. Then $\st
a \in Z(N)$ and $\st (\CS_a(R))= \CS_{\st a}(N)$.

\end{lemma}

\begin{proof}[Proof of the lemma]
Since $M\la a \ra\subseteq V$ the map $\st|_{M\la a \ra }:M\la a \ra
\to N$ is elementary. So $\st a \in Z(N)$.

Consider the function $f: Z\times M^n \to M$ given
 by $f(z,y) = d(\CS_z,y)$ (with $d(z,y) = \sqrt{z^2 + y^2}$
 and $d(S,y) = \inf \{d(z,y): z\in S\}$). So for each
 $z\in Z$ the function $f_z(y): f(z,y)$ is continuous.
  Hence by \S 6 (Section 2) of \cite{vdDr} there is a
  partition of $Z$ into definable subsets $Z_1,\dots,Z_k$
   such that each restriction $f_i|X_i\times M^n$ is continuous.
    Restricting $Z$, this allows us to assume that $f$ is continuous
     on $Z$. Let $b\in \CS_a(R)\cap V^n$; so $f(a,b)=0$. Since $f$ is
      continuous (and $\st a \in Z$), we get that $f(\st a, \st b) = \st f(a,b) = 0$.
       i.e. $\st b\in \cl \CS_{st a}$. Since we assume all the $\CS_z$ to be closed,
        this is what we wanted.

Conversely, if $c \in \CS_{\st a}$ we have that $\st c = c$ and therefore
 $f(\st a, \st c) =0$. Reversing the above argument, we get that $f(a,b)$
  is infinitesimal. By the definition of $f$ this means that $d(b,c)$
  (computed in $R$) is infinitesimal for some $b\in \CS_a(R)$. Thus,
   $c= \st b \in \CS_a(R)$, as desired.

\end{proof}

The rest of the proof is a rehash to what we did in Lemma
\ref{easycase}: By Proposition \ref{QE} the set $\tilde \CS
\setminus \CS$ is $\CM$-definable, and denote the corresponding
parameter set $Z'$.  It's enough to show that $\dim Z' < \dim Z$. We
take $\CM \prec \CN \prec \CR$ as in the previous lemma, and denote
$\st_M: R\to M$, $\st: R\to N$ the corresponding standard part maps.
As inthe proof of Lemma \ref{easycase} we obtain:
\[
 \{\tilde \CS_b: b\in Z'(M)\}=\st_M \CS_a\setminus \{\CS_a: a\in Z(M)\}
\]

Since $(\CR,\CM,\st_M) \equiv (\CR,\CN, \st)$ we also get
\[
 \{\tilde \CS_b: b\in Z'(N)\}=\st \CS_a\setminus \{\CS_a: a\in Z(N)\}
\]

By the previous lemma, if $b\in Z'(N)$ is generic and $a\in Z(R)$ is
such that
 $\st \CS_a(R) = \tilde \CS_b(N)$, it cannot be that $M\la a \ra \subseteq V$.
Thus, $\dim (\st M\la a \ra/M) < \dim (a/M)$, and since $b$ was
arbitrary, the desired result follows exactly as in the proof of
Lemma \ref{easycase}.

\section{Appendix B: The general o-minimal case}

In this appendix we adapt the main results of this paper to the
general o-minimal context (i.e. removing the assumption that $\CM$
expands a real closed field). From the proof of the results in
Section \ref{Tri} it is clear that we only need concern ourselves
with the proof of Theorem \ref{stable}.

In the first part of this appendix we will prove two important
technical facts for $\CN$ strongly minimal, 1-dimensional and
definable in a sufficiently saturated o-minimal structure $\CM$:
\begin{itemize}
 \item $\CN$ has finitely many non-orthogonality classes (in the o-minimal sense).
 \item $\CM$ has definable choice for $\CM$-definable subsets of $N^k$.
\end{itemize}

The second part of the appendix uses these two facts to overcome the
main obstacle in the generalisation of \ref{stable}: Lemma
\ref{asymptotes}. This lemma relies heavily on v.d. Dries' Theorem
on limit sets (\ref{limit sets}), which is only known to be true in
expansions of real closed fields.

Finally, we will go systematically through all places where the
existence of a field was used, and explain briefly how to avoid it.

Throughout this section $\CN$ will be a 1-dimensional structure {\em
definable} in a densely ordered  o-minimal $\CM$. It seems plausible
that a further generalisation of the theorem to the case where $\CN$
is only interpretable in $\CM$ can be achieved, but this will
definitely require more work.

\begin{claim} Without loss of generality, $N$ is a dense subset of
$M$.
\end{claim}
\proof  Since $\dim_{\CM} N =1$ it is a finite union of 1-cells and
0-cells (in $M^k$ for some $k$). Note that removing
finitely many points cannot alter the truth of Theorem
\ref{stable} and therefore neither of Theorem \ref{mainlocal} and Theorem \ref{main}.
 Hence we may assume that $N$ is
given by a union of 1-cells only. Each 1-cell is definably
homeomorphic to an interval in $M$ and therefore inherits an
ordering. Thus each 1-cell, with its induced structure from $\CM$
can be viewed as an o-minimal structure on its own. We now take
 the disjoint union of the 1-cells of $N$, ordering them together
 arbitrarily, so that each cell is an open interval in this
 ordering. We add endpoints between these intervals and obtain an
 o-minimal structure $\CM'$, such that $N$, without its
 0-cells, is a definable dense subset of $\CM'$.  We can now replace
$\CM$ by $\CM'$.\qed

\subsection{Non-orthogonality and definable choice}\label{defc}

 Since we can no longer assume that $\CM$ has elimination of
imaginaries, and as the work throughout the paper was carried out in
$\CM$ (and not in $\CM^{eq}$) we should be careful when using
$\CN$-interpretable subsets - as these will not, a priori, be
definable in $\CM$. Recall:

\begin{defnn}
A theory $T$ has \emph{weak elimination of imaginaries} if for every
$M\models T$ and $e\in M^{eq}$ there are $a_1,\dots, a_n\in M$ such
that $e\in \dcl(a_1,\dots,a_n)$ and $a_i\in \acl(e)$ for all $1\leq
i \leq n$.
\end{defnn}

It is fairly easy to check (See, e.g., \cite{HR2}) that:

\begin{lemma}
Assume $T$ is strongly minimal and $\acl(\0)$ is infinite then $T$
has weak elimination of imaginaries.
\end{lemma}

Hence, by adding constants to $\CN$ we may assume that it has weak
EI.

\bigskip

\begin{defnn}
Let $\CM$ be an o-minimal structure. {\em $a_1,a_2\in M$ are
non-orthogonal} if there exists a continuous monotone definable
bijection $f:I_{a_1} \to I_{a_2}$ for some open intervals
$I_{a_i}\ni a_i$ sending  $a_1$ to $a_2$ (see more in \cite{PeSt2}).

If $S\subseteq M$ is a definable set and $a\in M$ any element, we
say that {\em $b\not \perp a$ for every $b\in S$, uniformly in $b$},
if there exists a definable family of definable functions $F:M^2\to
M$ such that for all $x\in S$ the function $f_x(y):=F(x,y)$
witnesses that $x\not \perp a$.

\end{defnn}

The notion of (non) orthogonality plays an important role in
following arguments. The key observation, which makes everything
else work is:

\begin{lemma}\vlabel{orth}
If $\CN$ is non-locally modular then there exists $a\in \CN$ such
that $a\not \perp b$ uniformly in $b$ (in the o-minimal sense,
obviously) for all $b\in M$ outside a finite set.
\end{lemma}

\begin{proof}
Choose $a\in \CN$ generic (in the o-minimal sense). By weak EI we
can find $\CF$ an almost normal 2-dimensional family of
$\CN$-definable plane curves, as provided by non local modularity
(see Fact \ref{rich}). Hence for $b\in \CN$, o-minimally generic
over $a$, the point $\la a,b \ra$ is generic on a generic
$\CF$-curve through $\la a,b \ra$, implying that such a curve is,
locally near $\la a,b \ra$, the graph of a continuous monotone
function with $f(a)=b$, so $a\not \perp b$.

Consider now  the set $S$ of all $b\in \CN$ such that there is no
$\CF$-line through $\la a,b \ra$ which is locally near $\la a,b \ra$
the graph of a continuous monotone function. This is an o-minimally
definable set, and by what we have just shown, it must be
finite.\end{proof}

\begin{rmk}
The above lemma seems to remain true under the weaker assumption
that $\CN$ is non-trivial. In that case, if $\CN$ is locally
modular, use the fact that there exists an $\CN$-definable
equivalence relation $E$ with finite classes such that $N/E$ is a
co-finite set in a strongly minimal abelian group.
\end{rmk}

%


The following is easy:

\begin{lemma}\vlabel{triv}
If $\CM$ has finitely many non-orthogonality classes then $\CM$ does
not have a generic trivial type.
\end{lemma}

\begin{proof}
Assume that $e\in M$ is a generic trivial point. Let $I\ni e$
be a closed interval such that every point in $I$ is trivial (see,
e.g., the introduction to \cite{PeSt2}). We claim that every $a,b \in
I$ that are independent must be orthogonal to each other. Indeed, if
not then there is a definable, continuous, strictly monotone
function $f$ sending $a$ to $b$. It follows that for every $b'$ near
$b$ there is, uniformly in $b'$,  a definable $f_{ab'}$ witnessing
the non-orthogonality of $a$ and $b'$. It follows that $b$ is
nonorthogonal, uniformly, to every element in some neighbourhood of
it. It is not difficult to see that this contradicts the triviality
of $b$.\end{proof}

We can now prove that $\CM$ has definable choice:

\begin{lemma}\label{DefChoice}
 If $\CN$ is strongly minimal and non-trivial (and we still assume that $N$ is a
 dense subset of $M$) then $\CM$ has
definable Skolem functions. I.e. for every $\CM$-definable formula
  $\varphi(\x,\y)$ there is an $\CM$-definable function $f_{\varphi}$
   such that $\exists\y\varphi(\x,\y)\to \varphi(\x,f_{\varphi}(\x))$.
\end{lemma}

\begin{proof}
Since $M\setminus N$ is finite, it is enough to consider definable
subsets of $N^k$. By cell decomposition, we only need to consider
cells, and therefore, it suffices to prove that given a
$z$-definable interval $I_z$ we can find a $z$-definable point in
$I_z$.

If $\CN$ is locally modular, then since it is non-trivial there are
an $\CN$-definable vector space $V$, a finite equivalence relation
$E$ on $N$ and a definable injection $f: N/E \to V$ (see,
\cite{Zil5}) whose image is co-finite in $V$. By \cite{Ed},
$\CM$ has definable choice for definable subsets of $V$, this gives
an element $z_0\in f(I)$ and we take $z=\min\{f^{-1}(z_0)\}\in I$.

So we may assume that $\CN$ is not locally modular. Let
$(c,d)\subseteq N$ be any interval. We show that we can choose a
point $e\in (c,d)$ uniformly in $c,d$. Combining the previous lemma
with Lemma \ref{triv} for all but finitely $x\in N$ there exists,
uniformly in $x$, an interval $I_x$ on which a group-interval (not a
group!) is defined. If $c,d$ are not both of the finitely many
exceptional points, there is (without loss) an interval $I_c\ni c$
as above. If $d\notin cl(I_c)$ then setting $e:=\sup I_c\in (c,d)$
we are done. Otherwise $(c,d)\subseteq I_c$ and we can set $e=(c+_c
d)/_c 2_c$ (in the sense of the local group defined on $I_c$) if $d
< \sup I_c$ or $c+_c 1_c$ otherwise. Since the number of exceptional
points, and hence the number of intervals with exceptional
end-points, is finite, this gives us definable choice for sub
intervals in $N$.
\end{proof}

\subsection{Limit sets in structures without fields}

Our next task is to slightly refine the notion of non-orthogonality
in o-minimal structures.

\begin{defnn}
Let $\CM$ be an o-minimal (saturated enough) structure, $a,b\in
M\cup\{\pm\infty\}$. We say that $a^+\not \perp b^+$ if there exist
intervals $I_a^+:=(a,a_\epsilon)$, $I_b^+:=(b,b_{\epsilon})$ and a
definable (with parameters) continuous   monotone bijection $f:I_a^+
\twoheadrightarrow I_b^+$ (so $\lim_{t\to a}f(t)=b)$).
  For $\varepsilon_i\in \{+,-\}$ we define analogously,
  $a^{\varepsilon_1}\not \perp b^{\varepsilon_2}$ and $a^{\varepsilon}\not \perp \pm\infty$.
\end{defnn}

Consider the structure $\CM^{\pm}$ consisting of $\CM$ with all its
structure, as well as a new sort $M^{\pm}=(M\cup\{\pm\infty\})\times
\{+,-\}$ and the projection $\pi: M^{\pm}\to M\cup\{\pm \infty\}$.
 The relation $a^{\varepsilon_1}\not \perp b^{\varepsilon_2}$ is an
 equivalence relation
 on the elements of $M\times \{+,-\}$, so we can define:

\begin{defnn}
A complete set of representatives of the non-orthogonality classes
of $\CM^{\pm}$ is a subset of $\CM^{\pm}$ in which every
nonorthogonality class is represented.
\end{defnn}

\begin{defnn} Assume that $\CM$ is o-minimal and
every $a,b$ outside a finite set $F$ are non-orthogonal to each
other,  uniformly in $a,b$.  For the purposes of this note, we call
such $\CM$ \emph{uniformly unidimensional}.
\end{defnn}

\begin{lemma}
 Let $\CM$ be uniformly unidimensional, then $\CM$ has at most
  finitely many trivial points.
\end{lemma}
\begin{proof}
This is immediate from \ref{triv}.
\end{proof}

Note that by the previous lemma and (the proof of) Lemma \ref{DefChoice}
a uniformly unidimensional structure $\CM$ has definable choice. Moreover,
 uniformly unidimensional structures have curve selection for non-degenerate
  points. More precisely, if $S\subseteq M^r$ is a definable set and
  $q:=(q_1,\dots,q_r)\in \partial S$ is such that $q_i\notin F$, then
  there are $a < b\in M$ and a definable curve $\gamma:(a,b)\to S$ such
  that $\lim_{t\to a} \gamma(t) = q$.

\begin{defnn}
 Let $\CM$ be uniformly unidimensional. $a\in M$ is \emph{degenerate}  if
  one of $a^+, a^-$ is not in the generic non-orthogonality class of $\CM$.
\end{defnn}

\noindent A few clarifications may be in place:
\bigskip

\noindent{\bf Remark}

\noindent (1) If $\CM$ has finitely many non-orthogonality classes
so does $\CM^{\pm}$.  The converse is true as well: Let
$\{a_1^{\varepsilon_1},\dots, a_k^{\varepsilon_k}\}$ be a set of
representatives of the non-orthogonality classes of $\CM^{\pm}$. For
$b\in \CM$ define the non-orthogonality type of $b$ to be the pair
$(i,k)$ ($1\leq i,j\leq k)$ such that $b^+ \not \perp
a_i^{\varepsilon_i}$ and $b^- \not \perp a_j^{\varepsilon_j}$. It
will be enough to show that if $b,d$ have the same non-orthogonality
type, they are non-orthogonal in the usual sense, which is a simple
exercise in concatenation and composition of monotone functions.
\\

\noindent (2) The number of non-orthogonality classes, when finite,
depends only on $Th(M)$. Indeed, assume that
$\{a_1^{\varepsilon_1},\ldots, a_k^{\varepsilon_k}\}$ is a complete
set of representatives of the non-orthogonality classes of $\CM$.
Because of the saturation of $\CM$ the type stating that
$x^+\in M$ is orthogonal to each of the $a_i^{\varepsilon_i}$ is
inconsistent. It follows that, uniformly in $a$, every $a^+$ is
nonorthogonal to one of the $a_i^{\varepsilon_i}$, and similarly
for $a^{-}$. By quantifying over the $a_i$ we get in every model
$\CM_0$ of $Th(\CM)$ $k$-many non-orthogonality classes in
$\CM_0^{\pm}$.

\noindent (3) It is possible that $a^{\varepsilon_1}$ and
$b^{\varepsilon_2}$ are non-orthogonal to each other, while $a$ is a
degenerate point and $b$ is a non-degenerate point. For example,
consider $\la \mathbb R;<, +|(\mathbb R_{>0})^2\ra$.  We have
$0^+\not \perp r$ for every positive $r$ but $0$ is a degenerate
point.

\bigskip

\textbf{Fix $\CM$, a uniformly unidimensioinal structure}.

We now start our treatment of limit sets in $\CM$.
 Let $F^{\eps}=\{a_1^{\eps_1}, \ldots, a_k^{\eps_k}\}$ be the set of
 representatives of non-orthogonality classes in $\CM^{\pm}$
 which are orthogonal to $c^+$ for a fixed generic $c$. In particular,
 $F^{\eps}\cup\{c^+\}$ is a complete set of such representatives.
 We let $F=\{a_1, \ldots, a_k\}$. Clearly, $F\subseteq \dcl_{\CM}(\0)$.

 Recall that in the case of o-minimal expansions of fields,
 given a definable family of definable sets $\CS$, with parameter set $Z_1$,
  the limit family $\tilde \CS$ was defined as
   $\{\CS(\gamma): \gamma:(0,1)\to Z_1 \mbox{\, definable}\}$.
   In the presence of more than one non-orthogonality class not
   all definable curves $\gamma$ with image in $Z_1$ can have the same
   domain. In uniformly unidimensional structures, we handle
   separately limits along curves $\gamma:(a,b)\to Z_1$ with
   $a^+$ degenerate and those curves for which $a^+$ is non-degenerate.
Clearly, we can handle similarly curves $\gamma:(b,a)\to Z_1$ and
$a^-$.

Finally, it may happen that for some curve $\gamma$ the limit set
 $\CS(\gamma)$ is finite. To avoid unpleasant technicalities, and since
  in the application we will always be interested solely in limit sets that
  are infinite, we will restrict ourselves to those.
    We start with the former case:
\begin{lemma}\label{degenerate1}
 Let $\CM$ be be as above, $\CS$ an \emph{almost normal}
  family of definable curves with parameter set $Z_1$, and
   $\gamma:(a,b)\to Z_1$ a definable curve, $a^+\in F^{\eps}$.
   Assume that the limit set $\CS(\gamma)$ is infinite and  let
    $x=\la x_1,\dots,x_n\ra$ be a generic element in $\CS(\gamma)$.
    Then
     there exists some $i$ such that $x_i\in\dcl(\0)$. In particular,
    by genericity $x'_i = x_i$ for all $x'\in \CS(\gamma)\cap U$
     for some open $U\ni x$.
\end{lemma}

\begin{proof} As before we may assume that all the curves in $\CS$
are closed (after moving to the family of closures we still have an
almost normal family).
 Assume that the lemma fails. Then for some generic $x\in \CS(\gamma)$ all
  coordinates are non-orthogonal to some generic $c\in M$.
   Since $c$ is not trivial, there is, by the Trichotomy
    Theorem, a group interval defined around $c$. Pulling
     back this group structure to $x$ (on each coordinate
     separately) we may assume that $x$ lives in an expansion
      of a group interval. In particular, there is a definable family of open
       sets $U_r\ni x$, with $r\in M$ such that $\lim_{r\to c} U_r = \{x\}$ (the
       existence of this family is our only use of the group structure near $x$).
        Since $x\in \CS(\gamma)$, for every $r$ there exists $t(r)\in (a,b)$
        such that $S_{\gamma(t(r))}\cap U_r \neq \0$. Since $a^+\in F$, it
        cannot be that $t(r) \to  a$ as $r\to c$.  Therefore, $\lim_{r\to c}
        t(r) = a_0 > a$. It follows that $x\in S_{\gamma(a_0)}$.  But this is
        impossible: we can repeat the process, finding $a < ...<a_n < ...< a_1
         < a_0$ with $x\in S_{\gamma(t(a_i))}$ for all $i$, i.e. $x$ is in
         infinitely many of the $S_{\gamma(t)}$, so without loss, it appears
         in all of them (because this is a 1-dimensional family). Since we can
          repeat this for infinitely many $x$'s in $\CS_{\gamma}$ - this would contradict the
           normality of $\CS$.
\end{proof}

For $\gamma:(a,b)\to M^2$ a definable curve, we say that {\em
$\gamma$ is of rank $k$ over $A$} if the image of the interval
$(a,b)$ under $\gamma$ is a curve of rank $k$ over $A$.

From now one we also assume: {\bf Every generic type in $\CM$ is
rich, namely it is contained in an $\CM$-definable real closed
field}.

\begin{lemma}\vlabel{nondegenerate} Let $\CM$ be as above.
Let $\gamma=\la \gamma_1,\gamma_2\ra:(a,b)\to M^2$ be a definable
curve, with $a^+\notin F^{\eps}$. If $\gamma$ is of rank $k\geq 3$
over $A$ then there are open intervals $I_1, I_2\subseteq M$ such
that:
\\(i) For some $a<b_1<b$, the curve $\gamma(a,b_1)$ is contained in $I_1\times I_2$;
\\(ii) $I_1,I_2$ are definable over $B\supseteq A$ and the rank of
$\gamma|(a,b_1)$  over $B$ is still $k$.
\\(iii) $I_1, I_2$ are definably isomorphic to some open intervals
in a definable real closed field $R$.
\end{lemma}
\begin{proof} Since $\gamma$ has rank $k$ over $A$, there exists an
$A$-definable almost normal family of curves $\CX=\{X_d:d\in D\}$,
such that $\gamma(a,b)=X_d$ for some generic $d\in D$, and
$\dim(D)=k$. We first note that for every generic $x=\la
x_1,x_2\ra\in X_d$, we have $x_1, x_2\notin \dcl(\0)$. Indeed,
otherwise $X_d$ would be locally the graph of a constant function
and hence, in some $M^2$-neighbourhood of $x$, the family $\CX$ will
be 1-dimensional, contradicting the fact that $\CX$ is almost normal
and of rank at least $3$.

Let $\la p_1,p_2\ra\in (M\cup \{\pm \infty\})^2$ be the limit of
$\gamma(t)$ as $t$ tends to $a$. Because of our last remark, the
functions $\gamma_1, \gamma_2$ are nowhere constant hence for some
$\eps_1,\eps_2\in \{\pm\}$ we have $p_1^{\eps_1}=\lim_{t\to
a}\gamma_1(t)$ and $p_2^{\eps_2}=\lim_{t\to a}\gamma_2(t)$ (by
$p^{+} (p^{-})=\lim_{t\to a}\gamma_1(t)$ we mean that $\gamma(t)$
tends to $p$ from above (below)). To simplify notation, let us
assume that $\eps_1=\eps_2=+$. Since $a^+\notin F^{\eps}$, it is
nonorthogonal to $c^+$, for our fixed generic $c$. Let $R$ denote a
fixed definable real closed field such that $c\in R$.

There are now several cases to consider:

\noindent{\bf Case 1} $p_1, p_2\in \dcl(\0)$.

We consider intervals $I_1=(p_1,p_1')$, $I_2=(p_2, p_2')$, for
$p_1', p_2'$ generic and independent from $d$ over $A$. By choosing
$b_1$ sufficiently small, we have $\gamma(a,b_1)\subseteq I_1\times
I_2$. Because $\gamma_1,\gamma_2$ witness the fact that $p_1^+,
p_2^+\not \perp c^+$,  $I_1$, $I_2$ are definably isomorphic to open
intervals in the field $R$.

\noindent{\bf Case 2} $p_1,p_2\notin \dcl(\0)$.

In this case, $p_1$ and $p_2$ are nonorthogonal to $c$ and therefore
each is contained in a  definable real closed field. We can then
find open intervals $I_1\ni p_1, I_2\ni p_2$, definable over generic
parameters which are independent from $d$ over $A$, such that each
$I_i$ is contained in a definable real closed field. Because all
generic points are uniformly non-orthogonal to each other, the two
intervals are definably isomorphic to intervals in $R$.

\noindent{\bf Case 3} $p_1\in \dcl(\0)$ while $p_2\notin \dcl(\0)$.

In this case, we take $I_1=(p_1,p_1')$ and $I_2$ an interval
containing $p_2$ and proceed as before. We handle similarly the
symmetric fourth case.\end{proof}

\begin{lemma}\label{limitset} Let $\CM$ be as above.
 Let $\CS$ be a definable two-dimensional \emph{almost normal}
  family of curves, all contained in a two dimensional set $Q$,
  and having a parameter set $P\subseteq M^2$. Let
   $\gamma:(a,b)\to P$ be a definable curve of rank $k\geq 3$.

Assume that $S_{\gamma}$ is infinite.
   If $q$ is a generic element in $\CS(\gamma)$, with $\dim(q)=2$ then there
    exists an open neighbourhood $U$ of $q$ such that the rank of
$\CS(\gamma)\cap U$ is at most $2$.
\end{lemma}
\begin{proof} By Lemma \ref{degenerate1}, $a^+\notin F^{\eps}$, for otherwise
we will have $\dim(q)\leq 1$.

If $p=\lim_{t\to a}\gamma(t)$ then, by Lemma \ref{nondegenerate},
there exists $b_1\in (a,b)$ and an open rectangular box $V=
I_1\times I_2$ such that $\gamma(a,b_1)\subseteq V$, and $V$ is
definably isomorphic to an open rectangular box  in a definable real closed field $R$.
Moreover, the rank of $\gamma|(a,b_1)$ over all parameters is still
$k\geq 3$. Assume that $\gamma(a,b_1)=X_d$ for $X$ in a definable
family of curves $\CX=\{X_{d'}:d'\in D\}$. Then, after possibly
shrinking $D$ (but still with $\dim(D)=k$, we may assume that for
every $d'\in D$ the curve $X_{d'}$ is contained in $V$.

 Because $\dim(q)=2$, we may assume that $q=\la q_1,q_2\ra\in M^2$, where
 $q_1,q_2$ are nonorthogonal to $c$. Hence,
there exists an open rectangular box $U=J_1\times J_2\ni q$ such
that $J_1, J_2$ are definably isomorphic to open intervals in the
field $R$.

We can now restrict ourselves to $\CS(U,V)=\{\CS_{p'}\cap U:p'\in
V\}$. Because $U, V$ are definably isomorphic to open rectangular
boxes in the field $R$, we may assume that $\CS(U,V)$ and every
limit set of this family are  definable in $R$. We can now apply
Theorem \ref{limit sets} and conclude that the family of all possible limit
sets for the family $\CS(U,V)$ is definable of dimension at most
$2$. The limit curve $\CS(\gamma)\cap U$ belongs to this family and
therefore its rank is at most $2$.\end{proof}

Finally, we can prove Lemma \ref{asymptotes} without the field
assumption. We formulate it here again:
\begin{lemma} For $\CM$ uniformly unidimensional, let $\CF$ be an almost normal family
 of plane curves parameterised by $Q$, $X\subseteq P$ an $A$-definable curve
of rank $k>2$, $p_0\in X$ generic and assume that for $q_0\in Q$,
the curve $C_{q_0}$ touches $X$ at $p_0$ (here we identify $P$
locally with an open set in $M^2$), and in addition,
$\dim(p_0,q_0/\0)=3$. If $p$ is a special point for $(q_0,X)$ then
$\dim(p/A)=1$. In particular, $p$ is a generic point of $X$.
\end{lemma}

\begin{proof}
As was pointed out earlier on, $p\in \cl(X)$ and hence
$\dim(p/A)\leq 1$. By \ref{very normal}, we may assume that in a
neighbourhood $U$ of $p_0$ and $V$ of $q_0$, the family $\CF$ is
nice. Consider $\tau(X\cap U)$ and recall that $q_0$ is generic in
$\tau(X\cap U)$ over $A$ (by the Duality Theorem).

Assume, towards a contradiction, that $\dim(p/A)=0$. Hence, $q_0$ is
still generic in $\tau(X\cap U)$ over $Ap$. As was pointed out in
clause (3) of the discussion following Definition \ref{sp}, it
follows that $q_0$ is in the $Ap$-definable limit set $\CL(X,p)$.
But then, by the genericity of $q_0$ over $pA$, there exists a
neighbourhood $V_1\subseteq V$ of $q_0$, such that $\tau(X\cap
U)\cap V_1=\CL(X,p)\cap V_1$. In particular, $\CL(X,p)$ is infinite.

The limit set  $\CL(X,p)$ can also be written as $\CL(\gamma)$ with
$\gamma:(a,b)\to P$ a curve of rank $k>2$ over $A$ ($\gamma$
determines a branch of $X$ which contains $p$ in its frontier). By
Lemma \ref{degenerate1}, there exists an open neighbourhood $V'$ of
$q_0$, such that $\CL(\gamma)\cap V'$ has rank at most $2$ over $A$.
Without loss of generality, $V'=V_1$ and therefore the rank of
$\tau(X\cap U)\cap V_1$ over $A$ is at most 2. However, since the
rank of $X$ over $A$ is $k>2$, then the rank of $X\cap U$ is $k$ as
well and, by Lemma \ref{rank}, the rank of $\tau(X\cap U)\cap V_1$
is $k$ as well (see \ref{localrank}), contradiction.\end{proof}

\subsection{Some extra details}
In this final subsection, we explain - claim by claim - how to avoid
any use we made of the ambient field structure on $\CM$. As already
noted, the proof of Theorem \ref{1types} relies on results of the
first sections of this paper only in the treatment of the stable
case, and does not otherwise rely on $\CM$ expanding a field. It
will suffice therefore to show that Theorem \ref{stable} as well is
true independently of the field  assumption. Therefore, throughout
this section we will assume that $\CN$ is a strongly minimal
non-locally modular structure \emph{definable} in an arbitrary
o-minimal structure $\CM$.

The following lemma shows that under these (a posteriori
inconsistent) assumptions fields are (locally) definable in $\CN$:

\begin{lemma}
 Let $\CM$ be an o-minimal structure, $\CN$ a 1-dimensional non-locally
modular strongly minimal set definable in $\CM$ such that $N$ is
dense in $M$. Then $\CM$ is uniformly unidimensional and every
generic type in $\CM$ is rich.
\end{lemma}

\begin{proof} We already established that $\CM$ is uniformly unidimensional
(Lemma \ref{orth} and Lemma \ref{DefChoice}).

 Adding constants, we may assume that $\CN$ has weak EI. The
 non-local modularity of $\CN$ gives rise to a 2-dimensional family
 of plane curves in $\CN$ (and hence in $\CM$). Such a family cannot be
 definable in trivial structures or in ordered vector spaces, hence
 by the Trichotomy Theorem for o-minimal structures a real closed
 field $R$ is definable in $\CM$. By the unidimensionality of $\CM$
 every generic type in $\CM$ is nonorthogonal to a generic type in
 $R$.\end{proof}

From the last lemma it follows that all the results of Sections
\ref{Intersection} and \ref{SecStable} which uses the field
structure only in a local way transfer to the present context. Here
are some examples:
\begin{enumerate}
 \item In Lemma \ref{diff1}: The fact that $X$ is a graph of a
 $C^1$-function is meaningless without the existence of an ambient field. But
  since we are interested in the result only locally near a generic point
  $p_0=\la x_0, y_0 \ra \in X$ this can be solved as follows.
   Since $x_0, y_0\in N$
   are generic in $\CM$ they are non-orthogonal to each other and
   each has an $\CM$-definable real closed field containing it.
By non-orthogonality there exist closed (and non-trivial)
$\CM$-intervals $I_{1}\ni x_0, I_2\ni y_0$, which are definably
homeomorphic to each other. We may assume that $I_1, I_2$ are
definably homeomorphic to open intervals in the same field $R$.
 Hence, differentiation can be take with respect
       to that field structure. Note, also, that by stable embeddedness,
       $X\cap I_1\times I_2$ is definable, locally near $p_0$ using only
       parameters from $I_1\cup I_2$. Thus, the lemma clearly remains true
        in this new formulation.
 \item Lemma \ref{diff2}: Again, we use the fact that $x_0, y_0$ are
 nonorthogonal to generic elements in the same real closed field.
  \item The two lemmas discussed above are crucial to prove Theorem \ref{diff3},
  where the assumptions are only meaningful in the presence of an ambient field,
   but whose conclusion is independent of that field. Thus, we can reformulate
   the theorem as follows:
Let $\CF$ be a $\0$-definable (almost) normal family of plane curves
of dimension greater than $1$. Let $X$ be an $A$-definable plane
curve of rank greater than $1$, $p_0=\la x_0,y_0\ra\in X$ generic in
over $A$. Endow $\0$-definable neighbourhoods of $x_0$ and $y_0$
with a similar field structure, as above, and let $d$ be the derivative
- with respect to that field structure - of the function associated
to $X$ at $x_0$. If $p_0\in \Cq$ for $q_0\in Q$ and if
$f'_{q_0}(x_0)=d$, then $\Cq$ and $X$ touch each other at $p_0$ and
$\Mdim(p_0,q_0/\emptyset)=3$.
\end{enumerate}

This covers, all occurrences   of the field structure in Section
\ref{Intersection}. We now turn to Section \ref{SecStable}.

\begin{enumerate}
 \item The construction of the curve $X$, depends only on the $\CN$-structure
  and so causes no problem. We have to make sure, however, that
  there exists a curve from the family $\CF$ which touches
  $X$ at a generic point. This follows, essentially, from
 the version of Theorem \ref{diff3} we just formulated above. However, an
 alternative approach is also possible. Note that in the construction, we
 only need the points $z_1,z_2,x_0,y_0$ to be independent generics. By
 choosing them all very close to each other, and restricting our attention
 to a small interval containing them all, we may in fact assume
 (using stable embeddedness, as usual) that in that part of the
 proof we are working in an expansion of a field.
 \item Since we checked that the results of Section \ref{Intersection}
  go through, the proof of Theorem \ref{stable} in the first case, (the ``general case'')
   needs no alteration. As for the second part (the ``special case''), we
    note that the argument is only concerned with what is going on near
     (in the o-minimal sense) the point $p_0$. Thus choosing $z_1,z_2,x_0,y_0$
      as suggested above, will assure that the argument will need no
       alteration.
\end{enumerate}

\bibliographystyle{plain}

\bibliography{../../Bibfiles/harvard}

%
\end{document}